\documentclass[a4paper,reqno]{amsart}
\usepackage[utf8]{inputenc}
\usepackage[english]{babel}
\usepackage{amsmath,amssymb,amsthm,exscale,amsopn,mathtools,color,tikz-cd,cite,fullpage,enumitem,soul}
\usepackage[colorlinks=true,pdftex,unicode=true,linktocpage,bookmarksopen,hypertexnames=false]{hyperref}

\newtheorem{lem}{Lemma}[section]
\newtheorem{prop}[lem]{Proposition}
\newtheorem{thm}[lem]{Theorem}
\newtheorem{cor}[lem]{Corollary}

\newtheorem{theorem}{Theorem}
\newtheorem{proposit}[theorem]{Proposition}

\theoremstyle{definition}
\newtheorem{ex}[lem]{Example}
\newtheorem{rem}[lem]{Remark}

\newtheorem*{prob}{Problem}
\newtheorem*{assum}{Assumption}

\DeclareMathOperator{\Aut}{Aut}
\DeclareMathOperator{\ch}{char}
\DeclareMathOperator{\clK}{clKdim}
\DeclareMathOperator{\End}{End}
\DeclareMathOperator{\GK}{GKdim}
\DeclareMathOperator{\gld}{gldim}
\DeclareMathOperator{\gr}{gr}
\DeclareMathOperator{\id}{id}
\DeclareMathOperator{\Img}{Im}

\DeclareMathOperator{\lann}{lann}
\DeclareMathOperator{\M}{M}
\DeclareMathOperator{\Map}{Map}

\DeclareMathOperator{\pr}{pr}
\DeclareMathOperator{\Q}{Q}
\DeclareMathOperator{\rk}{rk}
\DeclareMathOperator{\Soc}{Soc}
\DeclareMathOperator{\Spec}{Spec}
\DeclareMathOperator{\Supp}{Supp}
\DeclareMathOperator{\Sym}{Sym}
\DeclareMathOperator{\Z}{Z}

\newcommand{\free}[1]{\langle #1 \rangle}

\newcommand{\mb}{\mathbb}
\newcommand{\mc}{\mathcal}

\newcommand{\ov}{\overline}
\newcommand{\s}{\subseteq}
\newcommand{\vn}{\varnothing}

\renewcommand{\le}{\leqslant}
\renewcommand{\ge}{\geqslant}

\newcommand{\tr}{\triangleright}

\author[I. Colazzo \and E. Jespers \and {\L}. Kubat \and A. Van Antwerpen]
{Ilaria  Colazzo \and Eric Jespers \and {\L}ukasz Kubat \and Arne Van Antwerpen}

\address[I. Colazzo (ORCID: 0000-0002-2713-0409)]{School of Mathematics, University of Leeds, Leeds, LS2 9JT, UK}
\email{I.Colazzo@leeds.ac.uk}

\address[E. Jespers (ORCID: 0000-0002-2695-7949)]
{Department of Mathematics and Data Science, Vrije Universiteit Brussel, Pleinlaan 2, 1050 Brussel}
\email{Eric.Jespers@vub.be}

\address[{\L}. Kubat (ORCID: 0000-0002-7848-6405)]{University of Warsaw, Institute of Mathematics, Banacha 2, 02-097 Warsaw, Poland}
\email{Lukasz.Kubat@mimuw.edu.pl}
\address[A. Van Antwerpen (ORCID: 0000-0001-7619-6298)]{Department of Mathematics and Statistics, National University of Ireland - Maynooth, Maynooth, Ireland}
\email{Arne.Vanantwerpen@mu.ie}

\title{Structure algebras of finite set-theoretic solutions of the Yang--Baxter equation}
\date{}

\begin{document}

\begin{abstract}
   Quadratic algebras related to some classes of finite  left non-degenerate solutions $(X,r)$ of the Yang--Baxter
   equation have been intensively studied since they are the associative ring-theoretical tool to study solutions. These
   are the monoid algebras $K[M(X,r)]$ and $K[A(X,r)]$, over a field $K$, of its structure monoid $M(X,r)$ and left derived
   structure monoid $A(X,r)$. In case $r$ is bijective (and thus also right non-degenerate) it is known that these algebras
   are representable (hence PI), left and right Noetherian and have finite Gelfand--Kirillov dimension. Moreover, such algebras
   are domains (or equivalently prime) if and only if they have finite global dimension, which also is equivalent to $r$ being
   an involutive map.\\
   In this paper we deal with structure algebras of arbitrary finite left non-degenerate solutions $(X,r)$, except for the last section.
   If $(X,r)$ satisfies additional conditions, such as being bijective, idempotent or left derived, it has been shown in a series of papers that
   $K[M(X,r)]$ is left Noetherian. In the first part of the paper we show that the algebra $K[M(X,r)]$ always is left Noetherian.
   Via divisibility by generators, we construct an ideal chain in $M(X,r)$ that has very strong algebraic
   structural properties on its Rees factors. This allows to obtain characterizations of when the algebras $K[M(X,r)]$ and $K[A(X,r)]$
   are right Noetherian. Intricate relationships between ring-theoretical and homological properties of these algebras and
   properties of the solution $(X,r)$ are proven. Furthermore, we describe the cancellative congruences of $A(X,r)$ and $M(X,r)$
   as well as the prime spectrum of $K[A(X,r)]$. This then leads to an explicit formula for the Gelfand--Kirillov dimension of
   $K[M(X,r)]$ and it equals the classical Krull dimension. 
    Finally, we obtain the first structural results for a class of finite degenerate solutions $(X,r)$ of the
   form $r(x,y)=(\lambda_x(y),\rho(y))$ by showing that structure algebras of such solution are always right Noetherian.
\end{abstract}

\maketitle

\section{Introduction}\label{sec:intro}

The study of the Yang--Baxter equation goes back to the papers of Yang \cite{Ya} and Baxter \cite{Ba} on statistical physics.
Let $V$ be a vector space over a field $K$. A solution of the Yang--Baxter equation is a linear map $R\colon V\otimes V\to V\otimes V$
satisfying the following equation on $V\otimes V\otimes V$:
\[(R\otimes{\id})({\id}\otimes R)(R\otimes{\id})=({\id}\otimes R)(R\otimes{\id})({\id}\otimes R).\]
Searching for solutions has attracted numerous studies both in mathematical physics and pure mathematics
and led to the introduction and investigation of  fundamental algebraic structures
including Hopf algebras and quantum groups, see for example \cite{Dri1992,Manin,KasselBook,BG}.

A description of all solutions of the Yang--Baxter equation is beyond current reach. For this reason, Drinfeld in \cite{Dri1992} suggested
focusing on those solutions that have an $R$-invariant basis $X$. In other words, solutions that are induced by a linear extension of a map
$r\colon X\times X \to X\times X$, where $X$ is a basis of~$V$, satisfying the following relation
\[(r\times{\id})({\id}\times r)(r\times{\id})=({\id}\times r)(r\times{\id})({\id}\times r)\] on $X\times X\times X$.
The pair $(X,r)$ is called a set-theoretic solution of the Yang--Baxter equation. Etingof, Schedler and Soloviev brought
algebraic tools into the study of such solutions \cite{ESS99}. Continuing in this direction, Rump introduced cycle sets
and braces \cite{Ru05,Rump2007} to provide an algebraic framework to generate and govern involutive non-degenerate
set-theoretic solutions, i.e. $r^2$ is the identity map. These tools then have been extended, by Guarnieri and Vendramin \cite{GV17} and Rump \cite{MR3881192},
in order to tackle general bijective non-degenerate solutions. These algebraic structures have since been widely studied, see for example 
\cite{MR3763276,MR4020748,MR3574204,MR3177933,MR3861714, MR4223285,MR3530867,MR2320986,Rump2007,MR3881192}, and there are strong and
fruitful connections with other areas, such as for example Hopf algebras \cite{MR3763907}, knot theory \cite{MR3868941}, group rings
\cite{ADS}, nil and nilpotent rings \cite{Sm}, and pre-Lie algebras \cite{MR4484785,bai2023post}.

The goal of this paper is to study quadratic algebras, which are structure algebras of set-theoretic solutions of the Yang--Baxter
equation, i.e., algebras over a field that are defined by quadratic relations determined by set-theoretic solutions of the
Yang--Baxter equation. Such algebras are called structure algebras and they are the ring-theoretical tool to study these solutions.
Interestingly, one of the original motivations to study quadratic algebras goes back to the theory of quantum groups, see for example
\cite{Manin} and \cite{MR2177131}. Structure algebras of set-theoretic solutions  have received a lot of interest in recent years,
mainly for finite non-degenerate bijective solutions, and this  in several contexts, such as arithmetical orders, classification of
four-dimensional non-commutative projective surfaces, Artin--Schelter regular algebras, regular languages and constructions of Noetherian
algebras (see for example \cite{CedoOkn2012,MR2927367,JO,JVC,YZ,MR3373388,COAut,MR2177131}). The first and motivating results tackling
quadratic algebras coming from set-theoretic solutions of the Yang--Baxter equation were obtained by Gateva-Ivanova and Van den Bergh
\cite{GVdB} for finite non-degenerate involutive  solutions. They showed that these algebras share many properties with polynomial
algebras in commuting variables, in particular they are representable (hence PI) algebras that are domains of finite global dimension,
which are left and right Noetherian and they also are subalgebras of group algebras of  solvable Bieberbach groups, i.e., a finitely
generated torsion-free abelian-by-finite groups. One of the aims of our work is to construct new classes of Noetherian algebras.
Because of their role in the algebraic approach in non-commutative geometry, low dimensional algebras (in terms of the homological
or the Gelfand--Kirillov dimension)  algebras  have gained a lot of attention and hence classes of algebras of this type are of
interest to researchers working in these areas. A second aim is to show that some properties of solutions are fully determined by
the algebraic properties of their structure algebras. 

We explain the above in greater detail. Let $X$ be a set and $r\colon X\times X\to X\times X$ a map.
Write \[r(x,y)=(\lambda_x(y),\rho_y(x)).\] The \emph{structure algebra},
over a field $K$, associated with such a map is, by definition,
\[\mc{A}_K(X,r)=K\free{X}/(xy-\lambda_x(y)\rho_y(x):x,y\in X),\] 
i.e., the quotient of the free $K$-algebra $K\free{X}$ on $X$ factored out by the ideal generated by the binomials
$xy-\lambda_x(y)\rho_y(x)$, with  $x,y\in X$. Since the presentation of the algebra is given by  homogeneous quadratic
word relations,  \[\mc{A}_K(X,r)\cong K[M(X,r)],\] i.e., the monoid algebra $K[M(X,r)]$, where $M(X,r)$ is the monoid
presented by generators from the set $X$ and subject to the defining relations $xy=\lambda_x(y)\rho_y(x)$ for all $x,y\in X$.
It is called the \emph{structure monoid} of $(X,r)$.

In the context of set-theoretic solutions $(X,r)$ of the Yang--Baxter equation, these monoids and algebras play a very crucial role and are being
intensely studied (see for example \cite{CJKVAV2020,CJV21,CJO2019,MR4278764,MR2927367,CoJeVAVe21x,JKVA2018,MR4105532,GatSeg,GVdB,gateva2011quantum}).
The strongest results have been obtained for finite non-degenerate bijective solutions. Recall that the solution $(X,r)$ is said to be finite
if $X$ is finite. It is called left non-degenerate (respectively right non-degenerate) if all maps $\lambda_x$ (respectively $\rho_x$) are bijective.
In case $(X,r)$ is both left and right non-degenerate one simply calls it a non-degenerate solution. Finally, $(X,r)$ is said to be bijective if $r$
is a bijective map.

We now outline the content of the paper. In Section~\ref{sec:background} we summarize the state-of-the-art and we give some preliminary results.
In Section~\ref{sec:lnoeth} and Section~\ref{sec:desc} we investigate the following problem for finite left non-degenerate solutions.

\begin{prob}
    When is the structure algebra $K[M(X,r)]$, over a field $K$, of a set-theoretic solution $(X,r)$
    of the Yang--Baxter equation left (respectively right) Noetherian?
\end{prob}

From now on we assume $X$ is finite as left or right Noetherianity of $K[M(X,r)]$ implies that $X$ is finite.
The first main result is the following: it settles the left Noetherianity problem, see Corollary~\ref{cor:thmA} in Section~\ref{sec:lnoeth}.

\begin{theorem}\label{theorem:A}
    Let $(X,r)$ be a finite left non-degenerate solution of the Yang--Baxter equation. If $K$ is a field, then $K[M(X,r)]$ is left Noetherian. 
\end{theorem}

In Section~\ref{sec:desc} we  tackle the problem of right Noetherianity for finite left non-degenerate solutions.
To do so we consider left divisibility in the monoid $M=M(X,r)$ by the generators $x\in X$ and we construct a finite ascending
ideal chain in $M$ of which the factors are either power nilpotent or  a uniform subsemigroup of a completely $0$-simple inverse
semigroup. Crucial for this construction are the sets $M_{YY}$ that consist of the elements $m\in M$ that are precisely left
divisible by the elements of a subset $Y$ of $X$ and satisfy $\lambda_m(Y)=Y$. As proven in Theorem~\ref{chainNoetherian},
Corollary~\ref{ArightNoetherian} and Proposition~\ref{rightnoethcancellative} these sets control when the algebra $K[M]$ is
right Noetherian. In Section~\ref{sec:lnoeth} we further relate the Noetherian problem of $K[M]$ to that of the algebra $K[S]$
for an easier, but closely related monoid $S=\Soc(M)=\{m\in M:\lambda_m=\id\}$, called the socle of $M$, and to the sets $S_{YY}=M_{YY}\cap S$.
Combining this with the results in the previous section one obtains (Proposition~\ref{MtoSoc} and Theorem~\ref{SrightNoetherian})

\begin{theorem}
    Let $(X,r)$ be a finite left non-degenerate solution of the Yang--Baxter equation. If $K$ is a field, $M=M(X,r)$ and $S=\Soc(M)$ then
    $K[M]$ is right Noetherian if and only if $K[S]$ is right Noetherian, which is further equivalent to the existence of a positive integer
    $d$ such that $S_{YY}^d$ is cancellative for each subset $Y$ of $X$ for which the Rees factor semigroup $(S_{YY}\cup M_{|Y|+1})/M_{|Y|+1}$
    is not nil.
\end{theorem}

In Section~\ref{sec:prime} we investigate primeness of $K[M]$ and it turns out that it is linked with a certain finite group
$\Omega_\lambda$, which is an orbit of a single element under the action of the permutation group $\free{\lambda_x:x\in X}$. We study
properties of $\Omega_\lambda$, we solve the problem of primeness of $K[M]$ when the group $\Omega_\lambda$ is trivial, and pose
the question whether this is always the case. Namely, summarizing Theorems~\ref{bijectiveqprime} and~\ref{primefinitegroup}, we obtain

\begin{theorem}
    Let $(X,r)$ be a finite left non-degenerate solution of the Yang--Baxter equation.
    If $K$ is a field and $M=M(X,r)$ then the following conditions are equivalent:
    \begin{enumerate}
        \item $(X,r)$ is an involutive solution.
        \item $M$ is a cancellative monoid and $\Omega_\lambda$ is a trivial group.
        \item $K[M]$ is a prime algebra and $\Omega_\lambda$ is a trivial group.
        \item $K[M]$ is a domain.
    \end{enumerate}
    Moreover, if the diagonal map $q\colon X\to X$, defined as $q(x)=\lambda_x^{-1}(x)$, is bijective then $\Omega_\lambda$ is a trivial group.
\end{theorem}

Let $A(X,r)=M(X,s)$ be the structure monoid of the left derived solution $(X,s)$ of $(X,r)$ given as $s(x,y)=(y,\sigma_y(x))$,
where $\sigma_y(x)=\lambda_y(\rho_{\lambda^{-1}_x(y)}(x))$. In Section~\ref{sec:semiprime} we characterize when the algebra
$K[A(X,r)]$ is semiprime in terms of a decomposition of the monoid $A(X,r)$ as a finite semilattice of certain cancellative
subsemigroups, see Theorem~\ref{semiprimeA} and Remark~\ref{remarkArch}. It turns out that $K[A(X,r)]$ is semiprime if and only
if all Archimedean components of $A(X,r)$ are cancellative. 

In Section~\ref{sec:cong} we first describe the cancellative congruence of $M(X,r)$ (Proposition~\ref{cancellativelnd}),
which is crucial in a description of prime ideals of $K[M(X,r)]$ having trivial intersection with $M(X,r)$. Next, we obtain
important information on all prime ideals of $K[A(X,r)]$ and this is then used to obtain a combinatorial formula for the
Gelfand--Kirillov dimension of $K[M(X,r)]$, see Theorem~\ref{GKdim}.

\begin{theorem}
    Assume $(X,r)$ is a finite left non-degenerate solution of the Yang--Baxter equation. Let $A=A(X,r)$ and $M=M(X,r)$.
    If $K$ is a field then \[\GK K[M]=\clK K[M]=\GK K[A]=\clK K[A].\] Moreover, for
    \[\mc{Z}_0=\{Z\s X:Z\ne X,\,\sigma_x(Z)\s Z\text{ and }\sigma_x(X\setminus Z)\s X\setminus Z\text{ for all }x\in X\setminus Z\},\]
    the above number is also equal to the maximum of numbers $t(Z)$ for $Z\in\mc{Z}_0$, where $t(Z)$ is the number of orbits
    of the set $X\setminus Z$ with respect to the action of the monoid $\Sigma_Z=\free{\sigma_x:x\in X\setminus Z}$.
    In particular, all the above (equal) dimensions are bounded by $|X|$.
 \end{theorem}

As a consequence of the previous result and the description of prime ideals in $K[A(X,r)]$ we obtain the following
homological characterization of involutive solutions in Corollary~\ref{cor:primefinitegroup}.

\begin{theorem}
    Let $(X,r)$ be a finite left non-degenerate solution of the Yang--Baxter equation.
    If $K$ is a field and $M=M(X,r)$ then the following conditions are equivalent:
    \begin{enumerate}
        \item $(X,r)$ is an involutive solution.
        \item $\GK K[M]=|X|$.
    \end{enumerate}
    Moreover, if the algebra $K[M]$ also is right Noetherian then the above conditions are equivalent to:
    \begin{enumerate}
        \setcounter{enumi}{2}
        \item $\rk M=|X|$.
        \item $\clK K[M]=|X|$.
        \item $\id K[M]=|X|$.
        \item $K[M]$ has finite global dimension.
        \item $K[M]$ is an Auslander--Gorenstein algebra with $\id K[M]=|X|$.
        \item $K[M]$ is an Auslander-regular algebra.
    \end{enumerate}
\end{theorem}

In Section~\ref{sec:allrhoequal}, we have a first result (Proposition~\ref{7.11}) on Noetherianity for some degenerate solutions.

\begin{proposit}
    Let $(X,r)$ be a finite solution of the Yang--Baxter equation of the form $r(x,y)=(\lambda_x(y),\rho(x))$.
    If $K$ is a field and $M=M(X,r)$ then $K[M]$ is a right Noetherian PI-algebra of finite Gelfand--Kirillov dimension.
    If, furthermore, the solution $(X,r)$ is left non-degenerate then $K[M]$ also is left Noetherian.
\end{proposit}

\section{Background and preliminary results}\label{sec:background}

We start this section with some general results known for set-theoretic solutions $(X,r)$ of the Yang--Baxter equation:
\begin{enumerate}[label=(R\arabic*),ref=R\arabic*]
    \item \label{R1} If $(X,r)$ is finite and left non-degenerate then $(X,r)$ is bijective if and only if $(X,r)$ is right non-degenerate
    (see \cite[Theorem~3.1]{CoJeVAVe21x}).
    \item \label{R2} If $(X,r) $ is finite and non-degenerate then the diagonal map $q\colon X\to X$, defined as $q(x)=\lambda_x^{-1}(x)$,
    is bijective (see \cite[Lemma~3.2]{CoJeVAVe21x}). Note that not all finite left non-degenerate solutions with a bijective
    diagonal map are bijective, as shown by the following example: $X=\{1,2\}$, $\lambda_1=\lambda_2=\id$ and $\rho_1(x)=\rho_2(x)=2$
    for each $x\in X$. 
    \item \label{R3} If $(X,r)$ is left non-degenerate then $(X,r)$ is bijective if and only
    if the left derived solution $(X,s)$ of $(X,r)$ is bijective or, equivalently, if the solution $(X,s)$ is right non-degenerate
    (see for example \cite[Proposition 2.2]{MR4105532}); actually $r = \varphi s\varphi^{-1}$ for some bijective map $\varphi$ on $X\times X$). 
\end{enumerate}

From now on, we adopt the following assumption, except for Section~\ref{sec:allrhoequal}.
\begin{assum}
    $(X,r)$ is a finite left non-degenerate set-theoretic solution of the Yang--Baxter equation.
    Moreover, we write $X=\{x_1,\dotsc,x_n\}$.
\end{assum}
Thus, $r\colon X\times X\to X\times X\colon(x,y)\mapsto(\lambda_x(y),\rho_y(x))$ is such that each
map $\lambda_x$ is bijective. Its associated structure monoid is
\begin{align*}
    M=M(X,r) & =\free{x_1,\dotsc,x_n\mid x_i \circ x_j=x_k\circ x_l\text{ if }r(x_i,x_j)=(x_k, x_l)}\\
    & =\free{x_1,\dotsc,x_n\mid x_i\circ x_j=\lambda_{x_i}(x_j)\circ\rho_{x_j}(x_i)}.
\end{align*}
The associated \emph{left derived solution} is \[s\colon X\times X\to X\times X\colon(x,y)\mapsto(y,\sigma_y(x)),\] 
where $\sigma_y(x)=\lambda_y(\rho_{\lambda^{-1}_x(y)}(x))$ and its associated structure monoid (called the \emph{left derived structure monoid})
is denoted $A=A(X,r)$.  Its operation is denoted additively. Let $\pi\colon M\to A$ be the bijective $1$-cocycle described in \cite[Proposition 3.2]{CJV21}
and put $a_i=\pi (x_i)$ for $1\le i\le n$. So \[A=A(X,r)=\free{a_1,\dotsc,a_n\mid a_i+a_j=a_j+\pi(\sigma_{x_j}(x_i))}.\]
Clearly, $A(X,r)$ and $M(X,r)$ are equipped with a natural length function, denoted $|\cdot|$, and the mapping $\pi$ is length
preserving. This length function also defines a gradation on the structure algebras $K[M(X,r)]$ and $K[A(X,r)]$.

It was shown by Gateva-Ivanova and Majid \cite[Theorem 3.6]{GIMa08} that for an arbitrary solution $(X,r)$ of the Yang--Baxter equation
the maps $\lambda\colon X\to\Map(X,X)$ and $\rho\colon X\to\Map(X,X)$ can be extended uniquely to $M=M(X,r)$ (for simplicity we use
the same notation for the extension) so that one  obtains a monoid homomorphism $\lambda\colon(M,\circ)\to\Map(M,M)\colon a\mapsto\lambda_a$,
and a monoid anti-homomorphism $\rho\colon(M,\circ)\to\Map(M,M)\colon a\mapsto\rho_a$. Moreover, the map $r_M\colon M\times M\to M\times M$,
defined as $r_M(a,b)=(\lambda_a(b),\rho_b(a))$, is a solution of the Yang--Baxter equation. Furthermore, for $a,b\in M$, we have
\begin{align*}
    a\circ b & =\lambda_a(b)\circ\rho_b(a),\\
    \rho_b(c\circ a) & = \rho_{\lambda_a(b)}(c)\circ\rho_b(a),\\
    \lambda_b(a\circ c) & =\lambda_b(a)\circ\lambda_{\rho_a(b)}(c).
\end{align*}
If the solution $(X,r)$ is left non-degenerate then it is proved in \cite[Proposition 2.11]{CoJeVAVe21x} that the structure
monoid $M=M(X,r)$ is a (unital) YB-semitruss $(M,+,\circ,\lambda,\sigma)$, with $\lambda$ defined as above, $\sigma\colon M\to\Map(M,M)$
given by $\sigma_b(a)=\lambda_b(\rho_{\lambda^{-1}_a(b)}(a))$, and $(M,+)$ isomorphic to $(A,+)$, where $A=A(X,r)$ is the left derived
monoid of $(X,r)$. Note that for any YB-semitruss $(B,+,\circ,\lambda,\sigma)$ it holds that $B+b\s b+B$, since
\begin{align}\label{sumsigma}
    a+b=b+\sigma_b(a)
\end{align}
for all $a,b\in B$. In particular, each right ideal of $(A,+)$ is a two-sided ideal of $A$. Moreover, $\lambda_a$ with $a\in M$
is an automorphism of $(M,+)$ and $\Img\lambda$ is a finite group. We refer to \cite{CoJeVAVe21x} for further details about YB-semitrusses.

The bijective $1$-cocycle $\pi\colon M\to A$ yields a monoid embedding \[f\colon M\to A\rtimes{\Img\lambda}\colon m\mapsto(\pi(m),\lambda_m).\]
Abusing notation, we often will identify $m$ with $f(m)$, i.e., we will write $m=(\pi(m),\lambda_m)$. So, in particular, we may also denote by
$x_i$ the generators of $A$.

So we simply may write \[M=\{(a,\lambda_a):a\in A\}=\free{x_i=(a_i,\lambda_{a_i}):1\le i \le n}.\]

We now  recall the following structural results known for the structure algebra $K[M(X,r)]$, over a field~$K$.

\begin{enumerate}[label=(S\arabic*),ref=S\arabic*]
    \item \label{S1} If $(X,r)$ is finite and left non-degenerate then $K[M(X,r)]$ is a representable algebra (i.e., it is embeddable
    in a matrix algebra over a field and thus PI). Indeed, we know by \cite[Theorem~5.3]{CoJeVAVe21x} that  $K[A(X,r)]$
    is a PI left Noetherian algebra. Thus, by a result of Ananin \cite{Ana1989}, it is representable. As $K[A(X,r)\rtimes{\Img\lambda}]$
    is a finitely generated free module over $K[A(X,r)]$, it is embeddable in the algebra of finite matrices over $K[A(X,r)]$.
    Hence the algebra $K[A(X,r)\rtimes{\Img\lambda}]$ is representable, and thus its subalgebra $K[M(X,r)]$ is representable as well.
    \item \label{S2} If $(X,r)$ is finite, bijective and non-degenerate then $K[M(X,r)]$ is a prime algebra if and only
    if it is a domain, and this is equivalent with $(X,r)$ being involutive (see \cite[Theorem~2.8]{JKVA2018} and \cite{MR4105532}).
\end{enumerate}

Of course, analogous definitions of structure monoid and derived structure monoid can also be given when $X$ is infinite.
However, it is worth noticing that a necessary condition for $K[M(X,r)]$ to be left or right Noetherian is that $X$ is finite,
because the commutative algebra $K[M(X,r)]/(xy:x,y\in X)$ is not Noetherian if $X$ is infinite. Moreover, certain restrictions
on $r$ are also required. Indeed, if $r=\id\colon X\times X\to X\times X$, then $K[M(X,r)]$ is the free $K$-algebra on the set
$X$ and thus is neither left nor right Noetherian if $X$ is not a singleton. 

The following key results provide
the state-of-the-art and our starting point on known properties of the structure algebra and derived structure algebra.

\begin{enumerate}[label=(N\arabic*),ref=N\arabic*]
    \item \label{N1} If $(X,r)$ is finite and left non-degenerate then $K[A(X,r)]$ is a left Noetherian algebra. Moreover,
    $\GK K[A(X,r)]\le |X|$ (see \cite[Theorem~5.3]{CoJeVAVe21x}). 
    \item  \label{N2} If $(X,r)$ is finite and left non-degenerate and, moreover, the diagonal map $q\colon X\to X$ is bijective
    then $K[M(X,r)]$ is left Noetherian and of finite Gelfand--Kirillov dimension  (see \cite[Corollary~5.4]{CoJeVAVe21x}).
    To prove this one shows that, under the assumptions, $K[M(X,r)]$ (and also $K[A(X,r)]$)
    is a finite left module over a left Noetherian algebra $K[B]$, for some submonoid $B$ of $A(X,r)$ such that $\lambda_b=\id$ for each
    $b\in B$, which ensures that the monoid $B$ may also be treated as a submonoid of $M(X,r)$.
    As an application of the result one immediately obtains that
    if additionally $(X,r)$ is bijective (and thus also right non-degenerate) then $K[M(X,r)]$ is a left and right Noetherian PI-algebra
    (see \cite{CoJeVAVe21x,JKVA2018,MR4105532}).
    \item \label{N3} If $(X,r)$ is finite, left non-degenerate and idempotent, that is $r^2=r$, then $K[M(X,r)]$ is left
    Noetherian and $K[M(X,r)]$ also is right Noetherian precisely when $|q(X)|=1$ (or equivalently, precisely when $M(X,r)$ is cancellative)
    (see \cite[Proposition~3.11 and Theorem~3.12]{MR4728712}).
\end{enumerate}

For a subset $S$ of $A$ we put, as in \cite{JKVA2018}, \[S^e=\{(s,\lambda_s):s\in S\}\s M.\]

We will use left divisibility to construct an ideal chain in $M=M(X,r)$. We follow the outline and construction
as in \cite[Section~3]{CJKVAV2020}. Recall that an element $s$ in a semigroup $T$ is said to be left divisible
by $t\in T$ if $s=t$ or $s=tt'$ for some $t'\in T$. We write $t\mid s$ in this case.

Note that, in $M$, an element $(a,\lambda_a)$ is left divisible by a generator $x_i=(a_i,\lambda_{a_i})$
if and only if $a$ is left divisible in $A$ by $a_i$. So, left divisibility in $M$ by elements of $X$ can be 
transferred to left divisibility in $A$ by elements of $\{a_1,\dotsc,a_n\}$, the generators of $A$.

For $1\le i \le n$ put \[M_i=\{m\in M:m\text{ is left divisible by at least $i$ different generators amongst }x_1,\dotsc,x_n\}.\]
Clearly, \[M_i=A_i^e=\{(a,\lambda_a):a\in A_i\},\] where
\[A_i=\{a\in A:a\text{ is left divisible by at least $i$ different generators amongst }a_1,\dotsc,a_n\}.\]
Note that $a\in A$ being left divisible by $a_i$ in $A$ means $a=a_i+b$ for some $b\in A$, or equivalently
\[(a,\lambda_a)=(a_i,\lambda_{a_i})\circ(\lambda_{a_i}^{-1}(b),\lambda_{a_i}^{-1}\lambda_a)
=x_i\circ(\lambda_{a_i}^{-1}(b),\lambda_{a_i}^{-1}\lambda_a).\]

Clearly, each $M_i$ is a right ideal of $M$. Because, by assumption, the solution $(X,r)$ is left non-degenerate and because,
for $a,b\in A$, we have $a+b=b+\sigma_b(a)$, it also is easily verified that each $M_i$ is a left ideal. Hence each $M_i$
is a two-sided ideal of $M$. Therefore, we get in $M$ the ideal chain
\begin{equation}
	\vn=M_{n+1}\s M_n\s M_{n-1}\s \dotsb \s M_1\s M_0=M.\label{Chain1}
\end{equation}
Next, we search for necessary and sufficient conditions in order to refine this chain to an ideal chain making
use of some intermediate ideals $B_i,U_i$ of $M$ satisfying
\begin{equation}
    M_{i+1}\s B_i \s U_i\s M_i\label{Chain2}
\end{equation}
and such that we have the following properties on the Rees factor semigroups $B_i/M_{i+1}$, $U_i/B_i$, $M_i/U_i$ and $M_i/M_{i+1}$:
\begin{enumerate}
	\item $B_i/M_{i+1}$ and $M_i/U_i$ are power nilpotent semigroups (if $M_i/M_{i+1}$ is power nilpotent then we take $B_i=U_i=M_i$),
	\item if $M_i/M_{i+1}$ is not power nilpotent then $U_i/ B_i$ is a disjoint union of semigroups $S_1,\dotsc,S_m$ such that 
	$S_k\circ S_l\s M_{i+1}$ for $k\ne l$, and 
	\item each $(S_i\cup M_{i+1})/M_{i+1}$ is a uniform subsemigroup of a completely $0$-simple inverse semigroup.
\end{enumerate}
Recall that a semigroup $T$ with a zero element $\theta$ is said to be power nilpotent if $T^k=\{\theta\}$ for some positive
integer $k$. Also recall that a completely $0$-simple inverse semigroup (also called a Brandt semigroup) is a semigroup of
the form $\mc{M}^0(C,k,k,I)$, where $C$ is a group and $I$ is the $k\times k$ identity matrix, i.e., this is the semigroup
of all $k\times k$ matrices with entries in $C^0=C\cup\{\theta\}$, the group $C$ with a zero $\theta$ adjoined, which have
at most one non-zero entry. A subsemigroup $T$ of $\mc{M}^0(C,k,k, I)$ is said to be uniform if each $\mc{H}$-class (i.e.,
all the matrices with non-zero entries in a fixed $(i,j)$ spot) of $\mc{M}^0(C,k,k,I)$ intersects non-trivially $T$ and
the maximal subgroups of $\mc{M}^0(C,k,k, I)$ are generated by their intersection with $T$. We make the agreement that
some ideals in the chain can be empty. For more details on semigroups we refer the reader to \cite{CP}.

The reason we would like to construct such an ideal chain is to determine when $K[M(X,r)]$ is right Noetherian by
making use of the following theorem of Okni\'nski \cite{MR1858039} (or see \cite[Theorem~5.3.7]{JO}).

\begin{thm}\label{chainNoetherian}
    Assume $S$ is a finitely generated monoid with an ideal chain \[\vn=S_{m+1}\s S_m\s S_{m-1}\s\dotsb\s S_1\s S_0=S\]
    such that each factor $S_j/S_{j+1}$ for $0\le j\le m$ is either power nilpotent or a uniform subsemigroup of a Brandt
    semigroup. Let $K$ be a field. If $S$ satisfies the ascending chain condition on right ideals and $\GK K[S]$
    is finite then $K[S]$ is right Noetherian.
\end{thm}

From \eqref{N1} we know that $\GK K[A]<\infty$, and thus (because $A=A(X,r)$ and $M=M(X,r)$ are in a bijection that
preserves the length function) also $\GK K[M]<\infty$, and that $K[A]$ is left Noetherian. In particular, the monoid
$A$ satisfies the ascending chain condition on left ideals. As each right ideal of $A$ is a two-sided ideal, we thus
also know that $A$ satisfies the ascending chain condition on right ideals. Now each right ideal of $M$ is of the
form $R^e$ for some right ideal $R$ of $A$. Hence $M$ satisfies the ascending chain condition on right ideals. 
So we can apply Theorem~\ref{chainNoetherian} once we have constructed a suitable ideal chain.
Because of its relevance we record these facts in the following lemma.

\begin{lem}\label{raccM}
    Assume $(X,r)$ is a finite left non-degenerate solution of the Yang--Baxter equation.
    Then both monoids $A(X,r)$ and $M(X,r)$ satisfy the ascending chain condition on right ideals.
    Furthermore, $A(X,r)$ satisfies the ascending chain condition on left ideals.
\end{lem}

Note that $M$ also satisfies the ascending chain condition on left ideals if, for example, $(X,r)$ also is right non-degenerate
(i.e., $(X,r)$ is bijective). Indeed for right non-degenerate solutions, one obtains the left-right version of the previous statement
by using the \emph{right derived monoid} (see \cite{CJV21})
\[A'(X,r)=\free{X\mid x+y=\rho_x(\lambda_{\rho^{-1}_y(x)}(y))+x\text{ for all }x,y\in X}.\]

\section{Left Noetherian algebras}\label{sec:lnoeth}

Let $(X,r)$ be a finite left non-degenerate solution of the Yang--Baxter equation. Let $A=A(X,r)$ and $M=M(X,r)$.
We begin by introducing the socle $\Soc(M)$ of $M$ as \[\Soc(M)=\{(a,\lambda_a):a\in A\text{ and }\lambda_a=\id\}.\]
Clearly this is a subsemigroup of $M$ and it is naturally isomorphic to a subsemigroup of $A$. Note that we do not
assume that the elements of the socle are central in $(A,+)$, contrary to what is done in the theory of skew braces
(see, e.g., \cite{GV17,MR4105532}). Consider the  diagonal map $q\colon X\to X$ with $ q(x)=\lambda_x^{-1}(x)$.
Since $X$ is finite, some power of the map $q$ is idempotent, that is there exist a positive integer $k$ such that 
\begin{equation}
    q^{2k}=q^k.\label{ekdef}
\end{equation}
For $x\in X$ and a positive integer $p$ we have \[(px,\lambda_{px})=(x,\lambda_x)\circ((p-1) q(x),\lambda_{(p-1)q(x)}).\]
Abusing notation, we will often write the generator $(q(x),\lambda_{q(x)})$ of $M$ simply as $q(x)$. It will be
clear from the context whether we consider this element in $M$ or in $A$. Then applying the previous equality
several times (for $p$ large enough) we get
\[(px,\lambda_{px})=x\circ q(x)\circ q^2(x)\circ\dotsb\circ q^{2k-1}(x)\circ q^k(x)\circ\dotsb\circ q^{2k-1}(x)\circ\dotsb\circ q^{p-1}(x).\]
In particular, \[(kq^k(x),\lambda_{kq^k(x)})=q^k(x)\circ\dotsb\circ q^{2k-1}(x).\] Hence for
\begin{equation}
    e=\exp\mc{G}(X,r),\label{expdef}
\end{equation}
the exponent of the permutation group $\mc{G}(X,r)=\free{\lambda_x:x\in X}$, we have
\[(keq^k(x),\lambda_{keq^k(x)})=(q^k(x)\circ\dotsb\circ q^{2k-1}(x))^e=(kq^k(x),\lambda_{kq^k(x)})^e.\]
Consequently, $\lambda_{keq^k(x)}=\id$. Thus \[(q^k(x)\circ\dotsb\circ q^{2k-1}(x))^e\in\Soc(M)\] and
\[(keq^k(x),\id)=q^k(x) \circ\dotsb\circ q^{k(e+1)-1}(x)=(q^k(x)\circ\dotsb\circ q^{2k-1}(x))^e.\] Put
\begin{equation}
    W=W_{ke}=\{(kex,\id):x\in X\text{ and }\lambda_{kex}=\id\}\s\Soc(M).\label{defW}
\end{equation}
and 
\begin{equation}
    Q=Q_{ke}=\{(keq^{k}(x),\id):x\in X\}\s W.\label{defQk}
\end{equation}
Note that the definition of $W$ depends on $ke$ but, for the results proven in the sequel, we may replace $ke$ by any of its multiples.

Let now $a\in A$ be an arbitrary element satisfying $|a|\ge nk(e+1)$ (recall that $n=|X|$). Then, for some $1\le i\le n$,
the generator $x_i$ appears at least $k(e+1)$ times in $a$, when $a$ is written as a sum of the generators $x_1,\dotsc,x_n$.
Hence \[a=k(e+1)x_i+b=b+\sigma_b(k(e+1)x_i)=b+k(e+1)x_j\] for some $b\in A$ and some $x_j\in X$. So, 
\[(a,\lambda_a)=(b,\lambda_b)\circ(\lambda_b^{-1}(k(e+1)x_j),\lambda_{\lambda_b^{-1}(k(e+1)x_j)})=(b,\lambda_b)\circ(k(e+1)x,\lambda_{k(e+1)x})\]
for some $x\in X$ and $|b|<|a|$. By the above 
\begin{align*}
    (k(e+1)x),\lambda_{k(e+1)x}) & =x\circ q(x)\circ\dotsb\circ q^{k-1}(x)\circ(q^k(x)\circ\dotsb\circ q^{2k-1}(x))^e\\
    & =x\circ q(x)\circ\dotsb\circ q^{k-1}(x) \circ w
\end{align*}
with $w=(q^k(x)\circ\dotsb\circ q^{2k-1}(x))^e=(kq^k(x),\lambda_{kq^k(x)})^e=(keq^k(x),\id)\in Q$. Hence
\[(a,\lambda_a)=(c,\lambda_c)\circ w\] for some $w\in Q$ and $c\in A$ with $|c|<|a|$.  Repeating the argument
on $c$ we thus  obtain  by induction  on the length of elements that \[M=\bigcup_{f\in F}f\circ \free{Q},\] 
where $F$ is the finite set of elements of $M$ of length less than $nk(e+1)$. In particular, $M$ is a finite right $\free{Q}$-module.

So we have shown the first part of the following proposition.

\begin{prop}\label{fromMtoS}
    Let $(X,r)$ be a finite left non-degenerate solution of the Yang--Baxter equation.
    If $K$ is a field, $M=M(X,r)$ and $S=\Soc(M)$ then the following properties hold:
    \begin{enumerate}
        \item $M=\bigcup_{f\in F}f\circ \free{Q}$, where $Q$ is defined in \eqref{defQk}.
        \item $S=\bigcup_{f\in F\cap S}f\circ \free{Q}=\bigcup_{f\in F\cap S} \free{Q}\circ f$,
        in particular $S$ is a finitely generated monoid.
        \item $M=\bigcup_{f\in F} \free{Q}\circ f$, and thus also $M=\bigcup_{f\in F} \langle W \rangle \circ f$.
        \item $(M,r_M)$ restricts to a solution on $W$, denoted by $(W,r_W)$. Furthermore, $(M,r_M)$ restricts
        to solutions on $S$ and $E=(F\cap S)\cup W$, denoted, respectively, by $(S,r_S)$ and $(E,r_E)$.
        \item $K[S]$ is a left Noetherian algebra and $S$ is an epimorphic image of the derived monoid $A(E,r_E)$.
        \item $S$ satisfies the ascending chain condition on left and right ideals.
        \item $K[\langle W \rangle ]$ is a left Noetherian algebra.
        Also $\langle W \rangle $ satisfies the ascending chain condition on left and right ideals.
    \end{enumerate}
    \begin{proof}
        To prove (2) fix $(s,\id)\in S$. Then, by part (1) of the Proposition, we have $(s,\id)=(a,\lambda_a)\circ(b,\id)$ for
        some $(a,\lambda_a)\in F$ and $(b,\id)\in \free{Q}$. But then $\lambda_a=\id$. Hence $(a,\lambda_a)\in F\cap S$, and
        thus $S=\bigcup_{f\in F\cap S}f\circ \free{Q}$. Moreover,
        \[(s,\id)=(a,\id)\circ(b,\id)=(a+b,\id)=(b+\sigma_b(a),\id)=(b,\id)\circ(\sigma_b(a),\id).\]
        Since $|\sigma_b(a)|=|a|$, we get $(\sigma_b(a),\id)\in F\cap S$, and thus $S=\bigcup_{f\in F\cap S} \free{Q}\circ f$.

        To prove part (3).
        Fix $(s,\lambda_s)\in M$. Then, by part  (1) of the Proposition, we have $(s,\lambda_s)=(a,\lambda_a)\circ(b,\id)$ for some
        $(a,\lambda_a)\in F$ and $(b,\id)\in \free{Q}$. But then $\lambda_a=\lambda_s$. Moreover,
        \[(s,\lambda_s)=(a,\lambda_s)\circ(b,\id)=(a+b,\lambda_s)=(b+\sigma_b(a),\lambda_s)=(b,\id)\circ(\sigma_b(a),\lambda_s).\]
        Since $|\sigma_b(a)|=|a|$, we get $(\sigma_b(a),\id)\in F$, and thus $M=\bigcup_{f\in F} \free{Q}\circ f$.
        
        To prove the first part of (4) note that for each $m=(key,\id)\in W$ and $w=(kex,\id)\in W$ we have
        $w\circ m=m\circ\rho_m(w)$. Moreover, writing $\rho_m(w)=(b,\lambda_b)$ for some $b\in A$ and comparing the
        $\lambda$-components in the previous equality we get $\id=\id\lambda_b$, and thus $\lambda_b=\id$. Furthermore,
        \[b=\lambda^{-1}_{key}(\sigma_{key}(kex))=\sigma_{key}(x))=ke(\sigma_{kex}(x))=kez,\]
        where $z=\sigma_{kex}(x)\in X$. Thus
        \begin{equation}\label{restrictw}
            \rho_m(w)=(b,\id)=(key,\id)\in W.
        \end{equation}
        The first part of (4) follows, because if $w,w'\in W$ then $\lambda_w(w')=w'\in W$ as well.
        
        To prove the second part of (4) choose $s=(a,\id)\in S$ and $t=(b,\id)\in S$. Then we have
        $r_M(s,t)=(t,\rho_t(s))=(t,\sigma_t(s))$. Next, note that
        \[(a+b,\id)=(b+\sigma_b(a),\id)=(b,\id)\circ(\sigma_b(a),\lambda_{\sigma_b(a)}).\]
        Hence $\lambda_{\sigma_b(a)}=\id$ and $\sigma_t(s)=(\sigma_b(a),\id)\in S$.
        Thus $r_M$ restricts to a solution on $S$. Furthermore, because of the first part of (3), equality 
        \eqref{restrictw} and the fact that each $\sigma_t$ respects the length, we obtain that $r_M$
        restricts to a solution on the set $E=(F\cap S)\cup W$.
        
        To prove part (5) we note that because of part (4), the monoid $S$ is an epimorphic image of the left derived
        monoid $A(E,r_E)$. Since the algebra $K[A(E,r_E)]$ is left Noetherian, also the algebra $K[S]$ is left Noetherian.
        In particular $S$ satisfies the ascending chain condition on left ideals. As right ideals of $S$ are two-sided
        ideals, $S$ also satisfies the ascending chain condition on right ideals, and thus (6) follows as well.

        To prove part (7) we note that because of part (4), the monoid $\langle W \rangle $ is an epimorphic image of the left derived
        monoid $A(W,r_W)$. Since the algebra $K[A(W,r_W)]$ is left Noetherian, also the algebra $K[\langle W \rangle ]$ is left Noetherian.
        In particular $\langle W \rangle $ satisfies the ascending chain condition on left ideals. As right ideals of $\langle W \rangle $
        are two-sided ideals, $S$ also satisfies the ascending chain condition on right ideals, and thus (6) follows as well.
    \end{proof}
\end{prop}

As an immediate consequence we obtain Theorem~\ref{theorem:A} (extending  results stated in \eqref{N1}, \eqref{N2} and \eqref{N3}).

\begin{cor}\label{cor:thmA}
    If $(X,r)$ is a finite left non-degenerate solution of the Yang--Baxter equation then the structure algebra $K[M(X,r)]$ is left Noetherian.
\end{cor}

The following example shows that the converse of the Corollary \ref{cor:thmA} does not hold.

\begin{ex}\label{ex3.3}
    Let $X$ be a finite set with at least two elements and choose a fixed element $z\in X$. Let $r\colon X^2\to X^2$ be the map
    defined by $r(x,y)=(z,z)$. This is a solution of the Yang--Baxter equation with all $\lambda_x =\rho_x$ equal to the constant
    map $c\colon X\to X$ given as $c(x)=z$. In $M(X,r)$ we have that every word of length $n\ge 2$ equals $z^n$. Hence,
    $K[M(X,r)]=K+\sum_{x\in X}Kx+\sum_{n\ge 2}Kz^n$ and thus $K[M(X,r)]$ is a finite left (and right) module over the commutative
    polynomial subalgebra $K[z]=\sum_{n\ge 0}Kz^n$. Since the latter is Noetherian, the algebra $K[M(X,r)]$ is left (and right)
    Noetherian, but $(X,r)$ is neither left nor right non-degenerate. Note that in this example all maps $\lambda_x$ and $\rho_x$
    are idempotent maps.
\end{ex}

Actually, Example \ref{ex3.3} is an illustration of a general result proven in Section~\ref{sec:allrhoequal}.
Indeed, the dual version of Proposition~\ref{7.11} says that if all $\lambda_x$ are equal, and not even necessarily
bijective, then $K[M(X,r)]$ is left Noetherian. 

\section{A description of the structure monoid and right Noetherian structure algebras}\label{sec:desc}

In this section we start by giving an explicit construction for an ideal chain for $M(X,r)$ (as in \eqref{Chain1} and \eqref{Chain2}).
For the construction we need a series of lemmas. Some of these are extensions of results that were proven in \cite{CJKVAV2020} for
bijective solutions, the proofs however, in  our  more general setting, become more tricky and involved. 

Fix $i$ with $1\le i \le n$. Let \[\mc{L}=\mc{L}(i)=\{Y\s\{a_1,\dotsc,a_n\}:|Y|=i\}.\] For $Y,Z\in\mc{L}$ put
\[M_{YZ}=\{(a,\lambda_a):a\in A\setminus A_{i+1},\,y\mid a\text{ for all }y\in Y\text{ and }\lambda_a(Z)=Y\}.\]
Notice that if $(a,\lambda_a)\in M_{YZ}$ then $a\in\free{Y}$ and $z\nmid a$ if $z\in \{a_1,\dotsc,a_n\}\setminus Y$.
Also note that some elements of $A_i \cap\free{Y}$ might belong to $A_{i+1}$. For $Y\in\mc{L}$ put
\[M_{Y*}=\bigcup_{Z\in\mc{L}}M_{YZ}\quad\text{and}\quad M_{*Y}=\bigcup_{Z\in\mc{L}}M_{ZY}.\]

\begin{lem}\label{lemma1}
	The following properties hold for $Y\in\mc{L}=\mc{L}(i)$:
	\begin{enumerate}
		\item $M_{Y*}\cup M_{i+1}$ is a right ideal of $M$.
		\item $M_{*Y}\cup M_{i+1}$ is a left ideal of $M$.
 	\end{enumerate}
 	\begin{proof}
		(1) Let $Z\in\mc{L}$ and $(a,\lambda_a)\in M_{YZ}$. Then, for $(b,\lambda_b)\in M$, we have
		$(a,\lambda_a)\circ(b,\lambda_b)=(a+\lambda_a (b),\lambda_a \lambda_b)$. If $a+\lambda_a(b)\notin\free{Y}$,
		i.e., $a+\lambda_a(b)$ is left divisible by some $z\in \{a_1,\dotsc,a_n\}\setminus Y$ then $a+\lambda_a(b)\in A_{i+1}$
		and thus $(a,\lambda_a)\circ(b,\lambda_b)\in M_{i+1}$. Otherwise, $a+\lambda_a(b)\in A_i\setminus A_{i+1}$.
		Furthermore, \[\lambda_{a+\lambda_a(b)}^{-1}(Y)=\lambda_b^{-1}(\lambda_a^{-1}(Y))=\lambda_b^{-1}(Z).\]
		Thus $(a,\lambda_a)\circ(b,\lambda_b)\in M_{Y\lambda_b^{-1}(Z)}\s M_{Y*}$, as desired.
		
		(2) Let $(a,\lambda_a)\in M_{ZY}$, where $Z\in\mc{L}$, and let $(b,\lambda_b)\in M$. Then
		$(b,\lambda_b)\circ(a,\lambda_a)=(b+\lambda_b(a), \lambda_b \lambda_a)$. If $b+\lambda_b(a)\in A_{i+1}$
		then $(b,\lambda_b)\circ(a,\lambda_a)\in M_{i+1}$. So, suppose $b+\lambda_b(a)\in A \setminus A_{i+1}$.
		Clearly, $\lambda_b(a)$ (and thus also $b+\lambda_b(a)$) is left divisible by all elements of $\lambda_b(Z)$.
		Thus $(b,\lambda_b)\circ(a,\lambda_a)\in M_i \setminus M_{i+1}$. Furthermore, we have
		\[(\lambda_b\lambda_a)^{-1}(\lambda_b(Z))=\lambda_a^{-1}(\lambda_b^{-1}(\lambda_b(Z)))=\lambda_a^{-1}(Z)=Y.\]
		Hence $(b,\lambda_b)\circ(a,\lambda_a)\in M_{\lambda_b(Z)Y}\s M_{*Y}$, as desired. 
	\end{proof}
\end{lem}

We fix the order of the elements of $\{a_1,\dotsc,a_n\}$ by their index. For $Y\in\mc{L}=\mc{L}(i)$ we thus obtain an
induced order on $Y$. So we write \[Y=\{a_{j_1},\dotsc,a_{j_i}\}\] with $1\le j_1<\dotsb< j_i\le n$, and we put
\[a_Y=a_{j_1}+\dotsb+a_{j_i}\in A.\] For $f\in\Sym(Y)$ we put \[a_{f,Y}=f(a_{j_1})+\dotsb+f(a_{j_i})\in A.\]
So $a_Y=a_{\id,Y}$. Now fix a positive integer $d$ (this number will be determined later) and put
 
\[m_{f,Y}=(da_{f,Y},\lambda_{da_{f,Y}})\in M\quad\text{and}\quad m_Y=m_{\id,Y}\in M.\]
 
Recall that a semigroup $T$ with a zero element $\theta$ is said to be nil provided each element of $T$
is nilpotent, that is for each $t\in T$ there exists a positive integer $k$ such that $t^k=\theta$. 
 
\begin{lem}\label{lemma2}
    Let $Y\in\mc{L}=\mc{L}(i)$. If $da_{f,Y}\in A_{i+1} $ (e.g., if $M_{YY}=\vn$) for some $f\in\Sym(Y)$
    and some positive integer $d$ then the following properties hold:
	\begin{enumerate}
		\item $(M_{YY}\cup M_{i+1})/M_{i+1}$ is nil.
		\item $(M_{*Y}\cup M_{i+1})/M_{i+1}$ is a nil left ideal of $M/M_{i+1}$.
		\item $(M_{Y*}\cup M_{i+1})/M_{i+1}$ is a nil right ideal of $M/M_{i+1}$.
	\end{enumerate}
	Hence if \[\mc{L}'=\{Y\in\mc{L}:da_{f,Y}\in A_{i+1}\text{ for some }f\in\Sym(Y)\text{ and }d\ge 1\}\]
    then $B_i=M_{i+1}\cup\bigcup_{Y\in\mc{L}'}(M_{*Y}\cup M_{Y*})$ is in the radical (i.e., the largest nil ideal) of $M/M_{i+1}$.
	\begin{proof}
        The statements to be proven are independent of the chosen order of the elements of $Y$. Hence for convenience, 
        it is sufficient to prove the lemma for $f=\id$.
		
		(1) Let $(a,\lambda_a)\in M_{YY}$. Then, for any positive integer $k$,
		\[(a,\lambda_a)^k=(a+\lambda_a(a)+\dotsb+\lambda_a^{k-1}(a),\lambda_a^k).\]
		Because $\lambda_a (Y)=Y$, we get that $(a,\lambda_a)^k\in M_{YY}\cup M_{i+1}$.
		Since each $\lambda_a^j(a)$ is left divisible by all elements of $Y$ and because for elements 
		$e,f\in A$, we have $e+f=f+\sigma_f (e)$,
		we get that $a+\lambda_a(a)+\dotsm+\lambda_a^{k-1}(a)$ is left divisible by $da_Y$ for
		a large enough $k$. Hence because of the assumption, it then follows that $(a,\lambda_a)^k\in m_Y\circ M\s M_{i+1}$. Therefore,
		$(M_{YY}\cup M_{i+1})/M_{i+1}$ is nil.
		
		(2) Let $Z\in\mc{L}$. Assume that $Z\ne Y$ and $(a,\lambda_a)\in M_{ZY}$. Then
		$(a,\lambda_a)\circ(a,\lambda_a)=(a+\lambda_a (a),\lambda_a^2)$. Because $\lambda_a(Y)=Z$
		and $Z\ne Y$ we have that $\lambda_a(Z)\ne Z$ and $a+\lambda_a (a)$ is left divisible by all
		elements in $Z\cup \lambda_a (Z)$. As $Z$ is properly contained in $Z\cup \lambda_a (Z)$ this
		yields $(a,\lambda_a)^2\in M_{i+1}$. Part (1) and Lemma~\ref{lemma1} therefore imply that
		$(M_{*Y}\cup M_{i+1})/M_{i+1}$ is a nil left ideal of $M/M_{i+1}$.
		
		(3) This is proved similarly as part (2).
	\end{proof}
\end{lem}

The following lemma deals with the case that $M_i/M_{i+1}$ is not nil.

\begin{lem}\label{lemma3}
    Suppose $Y\in\mc{L}=\mc{L}(i)$ and $(M_{YY}\cup M_{i+1})/M_{i+1}$ is not nil. Then:
    \begin{enumerate}
        \item if $d$ is a positive integer and $f\in\Sym(Y)$ then $da_{f,Y}\notin A_{i+1} $, i.e.,
        the generators of $A$ that left divide $da_{f,Y}$ are precisely those that belong to $Y$. 
        \item if $a_i,a_j\in Y$ are such that $a_i+a_j=a_k+a_l$ for some $a_k,a_l\in\{a_1,\dotsc,a_n\}$ then $a_k,a_l\in Y$.
        \item $A_{i+1}\cap \free{Y}=\vn$.
        \item for any $a\in A_i \setminus A_{i+1}$ divisible by all elements of $Y$, $f\in\Sym(Y)$ and positive integer $d$
        there exists $b\in A_i\setminus A_{i+1}$ such that $a+b=da_{f,Y}+b'$ and $b'\in A$ is a sum of elements of the form
        $da_{g,Y}$ with $g\in\Sym(Y)$.
        \item the left derived solution $s$ on $\{a_1,\dotsc,a_n\}$ restricts to a solution on $Y$, denoted by $s_Y$.
		\item $M_{YY}$ is a (non-empty) subsemigroup of $M$.
 	\end{enumerate}
	\begin{proof}
        (1) This follows at once from Lemma~\ref{lemma2}.
        
        (2) Put $x=a_i$, $y=a_j$, $u=a_k$ and $v=a_l$ and assume $x+y=u+v$. We need to show that $u,v\in Y$.
        From the assumptions we obtain that \[2a_Y=a'+x+y+b'=a'+u+v+b',\] for some $a',b'\in A$. Hence $2a_Y$
        is left divisible by $u$ and $v$. By part (1), $2a_Y\notin A_{i+1}$. Hence $2a_Y\in A_i$ and thus $u,v\in Y$, as desired.
        
        (3) This follows at once from part (2).
        
        (4) Let $a=b_1+\dotsb+b_t\in A_i\setminus A_{i+1}$ with each $b_j$ a generator of $A$ and each element of $Y$
        a divisor of $a$, and thus, by part (3), each $b_j\in Y$. Write $Y=\{a_{j_1},\dotsc,a_{j_i}\}$. As $a$ is left
        divisible by $f(a_{j_1})$, without loss of generality, we may assume that one of the summands of $a$ is $f(a_{j_1})$.
        Then we can rewrite $a=f(a_{j_1})+a'$ for some $a'\in A\setminus A_{i+1}$. Note that $|a'|<|a|$. It easily follows
        from parts (2) and (3)  that $a+a''=da_{f,Y}+a'''$ for some $a'',a'''\in A_i\setminus A_{i+1}$ with $|a'''|=|a'|<|a|$
        and the left divisors of $a'''$ belong to $Y$. The result now follows by repeating this argument on $a'''$.
		
		(5) This follows from part (2).
		
		(6) From Lemma~\ref{lemma2} we note that $M_{YY}$ is not empty. The statement now easily from part (2) and the definition of $M_{YY}$.
	 \end{proof}
\end{lem}

\begin{lem}\label{lemma4}
	Let $Y,Z,U,V\in\mc{L}=\mc{L}(i)$. If $Z\ne U$ then $M_{YZ}\circ M_{UV}\s M_{i+1}$.
	\begin{proof}
		Let $(a,\lambda_a)\in M_{YZ}$ and $(b,\lambda_b)\in M_{UV}$. Then $a+\lambda_a (b)$ is left divisible by all the elements
        of $ Y \cup \lambda_a (U)$. Because $Z\ne U$ and $\lambda_a (Z)=Y$, we have $\lambda_a (U) \ne Y$. Hence $|Y\cup\lambda_a (U)|>i$
        and thus $(a,\lambda_a)\circ(b,\lambda_b)=(a+\lambda_a (b),\lambda_a \lambda_b)\in M_{i+1}$, as desired.
	\end{proof}
\end{lem}

Let \[\mc{L}_u=\mc{L}_u(i)=\{Y\in\mc{L}:(M_{YY}\cup M_{i+1})/M_{i+1}\text{ is not nil}\}.\]
Because of Lemma~\ref{lemma3}, \[\mc{L}_u=\{Y\in\mc{L}:M_{YY}\text{ is a (non-empty) subsemigroup of }M\}.\]
We define a relation $\sim$ on $\mc{L}_u$ as follows. For $Y,Z\in\mc{L}_u$ we put
\[Y\sim Z\iff M_{YZ}\ne\vn \text{ or }M_{ZY}\ne\vn.\]

\begin{lem}\label{lemma6}
    Suppose $M_i/M_{i+1}$ is not nil (and thus $\mc{L}_u\ne \vn$) and $Y\in\mc{L}_u=\mc{L}_u(i)$.
    Then the following properties hold:
    \begin{enumerate}
        \item there exists $w_Y\in A_i\cap \free{Y}$ with $\lambda_{w_Y}=\id$. In particular, $(w_Y,\id)\in M_{YY}$.
        \item for each $a\in \free{Y}$ there exists $b\in \free{Y}$ and positive integer $d$ such that $a+b=dw_Y$,
        and clearly $\lambda_{a+b}=\id$.
    \end{enumerate}
    \begin{proof}
        (1) By assumption $M_{YY}$ is non-empty and it is a semigroup. Let $(a,\lambda_a)\in M_{YY}$ and let $k$ be the order of
        the permutation $\lambda_a$. Then $(a,\lambda_a)^k=(a+\lambda_a(a)+\dotsb+\lambda_a^{k-1}(a),\id)\in M_{YY}$.
        Hence it is enough to put $w_Y=a+\lambda_a (a)+\dotsb+\lambda_a^{k-1}(a)\in A_i\cap\free{Y}$.
        
        (2) Let $(w_Y,\id)\in M_{YY}$ and assume $a\in \free{Y}$ (and thus by Lemma~\ref{lemma3},
        $a$ is only left divisible by elements of $Y$). Write $a=a_j+c_j$ with $c_j\in\free{Y}$. Since $w_Y$
        is left divisible by each element of $Y$, we may write $w_Y=a_j+b_j$ for some $b_j\in\free{Y}$.
        Then, by Lemma~\ref{lemma3}, \[a+b_j=a_j+c_j+b_j=a_j+b_j+c'_j=w_Y+c'_j\] for some $c'_j\in\free{Y}$.
        Note that $|c'_j|<|a|$. So, repeating this argument on $c'_j$ we obtain that $a+b=dw_Y$ for some $b\in A$
        and some positive integer $d$ (as $|c_j'|=|a|-1$, one may take $d=|a|$).
    \end{proof}
\end{lem}

\begin{lem}\label{lemma5}
    Assume $M_i/M_{i+1} $ is not nil. Then the following properties hold:
	\begin{enumerate}
		\item if $Y,Z\in \mc{L}_u=\mc{L}_u(i)$ then $Y\sim Z$ if and only if $M_{YZ}\ne \vn$ and $M_{ZY} \ne \vn$.
		\item $\sim$ is an equivalence relation on $\mc{L}_u$.
	\end{enumerate}
	\begin{proof}
        (1) Let $Y,Z\in\mc{L}_u$ and $Y\sim Z$. Suppose $M_{YZ}\ne \vn$. Let $(a,\lambda_a)\in M_{YZ}$,
        in particular, $\lambda_a (Z)=Y$. By Lemma~\ref{lemma6} there exists $b'\in\free{Y}$ and a positive integer
        $d$ so that $a+b'=dw_Y$. Hence if $b=\lambda_a^{-1}(b')$ then
        \[(a+\lambda_a(b),\lambda_a\lambda_b)=(a,\lambda_a)\circ(b,\lambda_b)=(dw_Y,\id)\in M_{YY}.\] 
		If necessary, multiplying with another factor $(w_Y,\id)$, we may assume that
        $b\in A_i\setminus A_{i+1}$. As $\lambda_a \lambda_b=\id$ we have that $\lambda_b=\lambda_a^{-1}$. In particular,
        $\lambda_b (Y)=Z$. Because of Lemma~\ref{lemma3}, the generators of $A$ that left divide $\lambda_a (b)$ are precisely
        the elements of $Y$. Hence the generators of $A$ that left divide $b$ are precisely the elements of $\lambda_a^{-1}(Y)=Z$.
        It follows that $(b,\lambda_b)=(b,\lambda_a^{-1})\in M_{ZY}$. Hence $M_{ZY}\ne \vn$. Part (1) then follows.
		
		(2) Clearly $\sim$ is reflexive and symmetric. To show that it is transitive, let $Y,Z,U\in\mc{L}_u$ with
		$Y\sim Z$ and $Z\sim U$. So, by (1), the sets $M_{YZ},M_{ZY},M_{ZU},M_{UZ}$ are all non-empty. It follows that also all
		of the sets $(w_Y,\id)\circ M_{YZ}\s M_{YZ}$, $(w_Z,\id)\circ M_{ZY}\s M_{ZY}$, $(w_Z,\id)\circ M_{ZU}\s M_{ZU}$,
        $(w_U,\id)\circ M_{UZ}\s M_{UZ}$ are non-empty. We need to show that the set $M_{YU}$ is non-empty. To do so,
        it is sufficient to show that ($(w_Y,\id)\circ M_{YZ})\circ((w_Z,\id)\circ M_{ZU})$ is not contained in $M_{i+1}$.
        Let $(a,\lambda_a)\in M_{YZ} $ and $(b,\lambda_b)\in M_{ZU} $. Since $\lambda_a(Z)=Y$ then $\lambda_a(w_Z)$ is
        precisely left divisible by all elements of $Y$. Note that
        \begin{align*}
            (w_Y,\id)\circ(a,\lambda_a)\circ(w_Z,\id)\circ(b,\lambda_b) & =(w_Y+a,\lambda_a)\circ(w_Z+b,\lambda_b)\\
            & =(w_Y+a+\lambda_a (w_Z)+\lambda_a(b),\lambda_a\lambda_b).    
        \end{align*}
		Because of Lemma~\ref{lemma3}, the generators of $A$ that left divide $w_Y+a+ \lambda_a (w_Z)+\lambda_a (b)$
        are precisely the elements of $Y$. Furthermore, $\lambda_a \lambda_b (U)=\lambda_a (Z)=Y$. Hence
		$(w_Y,\id)\circ(a,\lambda_a)\circ(w_Z,\id)\circ(b,\lambda_b)\in M_{YU}$, as claimed. This proves part (2).
	\end{proof}
\end{lem}

To prove the next lemma we will make use of \cite[Theorem~3.5]{Okn98} on the structure of linear monoids, i.e.,
multiplicative subsemigroups, say $T$, of a full linear monoid of square matrices $\M_m(L)$ over a field $L$. 
The matrices in $T$ of rank at most $j$ are denoted by $T_{(j)}$. It is shown, in particular, that such a monoid $T$
has a finite ideal chain with each Rees factor either a power nilpotent semigroup or a $0$-disjoint union of uniform
subsemigroups of a completely $0$-simple semigroup. In case the semigroup $T$ does not have a zero element, it follows
that the lowest ideal in the chain (that is the elements of $T$ of minimal rank) is a uniform subsemigroup of
a completely simple semigroup.

\begin{lem}\label{lemma8}
    Suppose $M_i/M_{i+1}$ is not nil and $Y\in\mc{L}_u=\mc{L}_u(i)$. Let $(w_Y,\id)\in M_{YY}$ be as in Lemma~\ref{lemma6}.
    Then there exist positive integers $d$ and $k$ such that $M_{YY}^k\s (dw_Y,\id)\circ M_{YY}$ is a subsemigroup of a
    completely simple semigroup $\mc{M}(G,1,V,P)$. Thus $(dw_Y,\id)\circ M_{YY}=\bigcup_{v\in V}S_{Y,v}$ is a disjoint
    union of a family of cancellative semigroups $(S_{Y,v})_{v\in V}$ such that $S_{Y,v}\circ S_{Y,v'}\s S_{Y,v'}$
    for all $v,v'\in V$. In particular, each $S_{Y,v}$ satisfies the left and right Ore condition, and $S_{Y,v}$ has
    a group of fractions $G_{Y,v}$ that is abelian-by-finite and isomorphic with $G$.
    Consequently, $(dw_Y,\id)\circ M_{YY}$ is left cancellative.
    
    Furthermore, $(dw_Y,\id)\circ M_{YY}$ is cancellative if and only if there exists a positive integer $N$ so that
    $a^N \circ b^N=b^N\circ a^N$ for all $a,b\in (dw_Y,\id)\circ M_{YY}$, or equivalently for any $a\in (dw_Y,\id)\circ M_{YY}$ 
    one has $((dw_Y,\id)\circ M_{YY})\cap((dw_Y,\id)\circ M_{YY}\circ(dw_Y,\id))\ne\vn$ (it actually is sufficient that this condition
    holds in $A$).
    \begin{proof}
        If $K$ is a field then $K[M]$ is a representable algebra by \eqref{S1} and thus $M$ is a linear monoid.
        By Lemma~\ref{lemma3}(6) we know that $M_{YY}$ is a subsemigroup of $M$, and it clearly does not have
        a zero element (as it does not contain idempotent elements of length at least one). Hence $M_{YY}$ is
        a linear semigroup without a zero element, and consequently, we may consider $M_{YY}$ as a subsemigroup
        of the multiplicative monoid $\M_m(L)$, for some field $L$. Let $S$ be the set of elements of $M_{YY}$ that
        have minimal rank (considered as elements of $\M_m(L)$). By the remarks stated before the lemma, we know that $S$ is
        a uniform subsemigroup of a completely simple semigroup. In particular, $S$ is a disjoint union of cancellative semigroups.
        Let $(w_Y,\id)\in M_{YY}$ denote the element described in Lemma~\ref{lemma6}. In particular, $(w_Y,\id)\circ M_{YY}\circ S\s S$.
        Let $m'\in (w_Y,\id)\circ M_{YY}\circ S$. Because of Lemma~\ref{lemma6}(2), there exist a positive integer $d$ and $m''\in M_{YY}$
        so that $m_Y=m'\circ m''=(dw_Y,\id)\in M_{YY}\circ S$. Hence $m_Y\circ M_{YY}$ is a subsemigroup of $S$ and thus a disjoint union
        of cancellative semigroups, say $m_Y\circ M_{YY}=\bigcup_{k,l}S_{kl}$, where each $S_{kl}$ is a subsemigroup of a maximal
        subgroup of the completely simple semigroup in which $S$ is uniform and each maximal subgroup of the completely simple
        semigroup is generated by its intersection with $m_Y\circ M_{YY}$, i.e., by one of the semigroups $S_{kl}$. Since the
        algebra $K[S_{kl}]$ is PI we know (see, e.g., \cite[Theorem~3.1.9]{JO}) that $S_{kl}$ has a group
        of fractions that is abelian-by-finite. In particular, $S_{kl}$ is a left and right Ore semigroup. Note that
        $S_{kl}\circ S_{k'l'}\s S_{kl'}$ for all $k,l,k',l'$. Again by a previous comment, every element of $m_Y\circ M_{YY}$
        has a right multiple of the type $m_Y^j=(jdw_Y,\id)$ (that belongs to $M_{YY}$). Hence all elements of $m_Y\circ M_{YY}$
        have the same row index as the element $m_Y$. So there only is one first index, say $1$, and we denote $S_{1l}$ as $S_{Y,l}$.
        Since elements of large enough length in $M_{YY}$ belong to $m_Y\circ M_{YY}$, the first part of the result follows.
        
        For the second part, we note that $(dw_Y,\id)\circ M_{YY}=\bigcup_lS_{Y,l}$ is right cancellative if and only if there
        only is one index $l$. The claim now follows from the previous comments and the fact that the group of quotients of
        the cancellative component containing $m_Y$ is obtained by inverting $m_Y$.
    \end{proof}
\end{lem}

It is easy to give examples of finite left non-degenerate solutions $(X,r)$ of the Yang--Baxter equation for which not
all $(dw_Y,\id)\circ M_{YY}$ are cancellative.

\begin{ex}\label{ex:idem}
    Assume $(X,r)$ is a finite left non-degenerate idempotent solution of the Yang--Baxter equation with $|X|>1$. Then in $A=A(X,r)$
    one has $a+b=b+b$ for all $a,b\in A$ of equal length. Hence any element of length at least two is left divisible by all elements of $X$.
    So, \[A_{XX}=A\setminus (X\cup\{1\})=\bigcup_{x\in X}(\free{x}+2x)\quad\text{and}\quad dA_{XX}=\bigcup_{x\in X}d(\free{x}+2x).\]
    This semigroup clearly is not right cancellative.
\end{ex}

We now state the second main result. 

\begin{thm} \label{idealchainNoetherian}
    Let $(X,r)$ be a finite left non-degenerate solution of the Yang--Baxter equation. Let $K$ be a field and $M=M(X,r)$. 
    For each subset $Y$ of $X$ for which $M_{YY}$ is a (non-empty) semigroup, i.e., $Y\in\mc{L}_u(|Y|)$, choose $(w_Y,\id)\in M_{YY}$
    as in Lemma~\ref{lemma6}. If for each such subset $Y$ of $X$ there exists a positive integer $d$ so that $(dw_Y,\id)\circ M_{YY}$
    is a cancellative semigroup then $K[M]$ is right Noetherian. The cancellative assumption holds for a (non-empty) semigroup
    $(dw_Y,\id)\circ M_{YY}$ if and only if for some positive integer $N$ one has $a^N\circ b^N=b^N\circ a^N$ for all
    $a,b\in (dw_Y,\id)\circ M_{YY}$.
    \begin{proof}
        The second  part of the statement follows from Lemma~\ref{lemma8}.
        
        To prove the first part, assume $(dw_Y,\id)\circ M_{YY}$ is a cancellative semigroup for each $Y$ for
        which $M_{YY}$ is a (non-empty) semigroup. The result follows from all other previous lemmas (as was
        done in the proof of \cite[Proposition~3.6]{CJKVAV2020} for bijective left non-degenerate solutions)
        and the comments stated after Theorem~\ref{chainNoetherian}. For completeness' sake and clarity to the
        reader we explicitly state the chain of ideals that satisfies the properties stated in \eqref{Chain1}
        and \eqref{Chain2}. By Theorem~\ref{chainNoetherian} and because $M$ satisfies the ascending
        chain condition on right ideals by Lemma~\ref{raccM}, it follows that $K[M]$ is right Noetherian. If $M$ also satisfies
        the ascending chain condition on left ideals then it follows that $K[M]$ also is left Noetherian. 
        
        Fix $1\le i \le n$ (recall that $n=|X|$). Let $\mc{L}_1 ,\dotsc,\mc{L}_{k}$ denote the equivalence classes of
        the relation $\sim$ on $\mc{L}_u=\mc{L}_u(i)$ (see Lemma~\ref{lemma5}). For each $1\le j \le k$ put
        \[\mc{U}_{ij}=\bigcup_{Y,Z\in\mc{L}_j}M_{YZ},\quad U_{ij}=\bigcup_{Y,Z\in\mc{L}_j}m_Y\circ M_{YZ},\quad U_i=\bigcup_{j=1}^kU_{ij}.\]
        The following properties hold:
        \begin{enumerate}
            \item $(\mc{U}_{ij} \cup M_{i+1})/M_{i+1} $ is a subsemigroup of $M_i/M_{i+1}$ such that
            $M_{YZ}\circ M_{ZV}\s M_{YV}$ and $M_{YZ}\circ M_{UV}\s M_{i+1}$ for all $Y,Z,U,V\in\mc{L}_j$ with $U\ne Z$.
            \item $(U_{ij} \cup M_{i+1})/M_{i+1}$ is an ideal of $(\mc{U}_{ij} \cup M_{i+1})/M_{i+1} $ and it is
            a subsemigroup of a completely $0$-simple inverse semigroup with maximal subgroups the group of fractions
            of $(dw_Y,\id)\circ M_{YY}$.
            \item $(\mc{U}_{ij} \cup M_{i+1})/M_{i+1} $ does not contain a nil ideal.
            \item $B_i=M_{i+1}\cup\bigcup_{Y\in\mc{L}\setminus \mc{L}_u}(M_{*Y}\cup M_{Y*})$ is the radical of $M_i/M_{i+1}$.
            Furthermore, if $Y\in\mc{L}\setminus\mc{L}_u$ and $Z\in\mc{L}_u$ then $M_{YZ}=\vn$ or $M_{ZY}=\vn$.
            \item $M_i/B_i=\bigcup_{j=1}^k(\mc{U}_{ij} \cup B_i)/B_i$, a $0$-disjoint union.
            \item $M_i/(U_i \cup M_{i+1})$ is a nil semigroup.
        \end{enumerate}
        Hence in $M$ we have an ideal chain
        \[M_{i+1}\s B_i\s U_{i1}\cup B_i\s U_{i1}\cup U_{i2}\cup B_i\s
        \dotsb\s U_{i1}\cup U_{i2} \cup \dotsb \cup U_{ik}\cup B_i=U_i\cup B_i\s M_i,\]
        where the first and last Rees factor is a power nilpotent semigroup and all other Rees factors are uniform
        subsemigroups of an inverse completely $0$-simple semigroup with maximal subgroups the groups of fractions
        of cancellative subsemigroups of $M$.
    \end{proof}
\end{thm}

\begin{cor}\label{ApplOre}
    Let $(X,r)$ be a finite left non-degenerate solution of the Yang--Baxter equation. Let $K$ be a field,
    $A=A(X,r)$ and $M=M(X,r)$. If $A$ satisfies the left Ore condition (e.g., if $A$ is abelian, or more generally
    a Malcev nilpotent monoid, see \cite{JO,CJKVAV2020} for the definition) then $K[M]$ is right Noetherian.
    \begin{proof}
        This follows at once from Theorem~\ref{idealchainNoetherian} because the assumptions easily imply that
        for $Y\in\mc{L}_u=\mc{L}_u(i)$ the semigroup $ (dw_Y,\id)\circ M_{YY}$ is cancellative by Lemma~\ref{lemma8}.
    \end{proof}
\end{cor}

Recall that in \cite{MR2290908} regular submonoids of the holomorph of a cancellative abelian monoid
have been studied, i.e., submonoids $M=\{(a,\lambda_a):a\in A\}$ of $A\rtimes{\Aut(A,+)}$ with $A$ a cancellative
abelian monoid. Such monoids are  said to be of IG-type \cite{MR2290908}; in case $A$ is free abelian they were
studied earlier by Gateva-Ivanova and Van den Bergh \cite{GVdB} and Jespers and Okni\'nski \cite{MR2189580} and
are called monoids of I-type. In the latter case, their monoid algebras are the structure algebras of involutive
non-degenerate solutions of the Yang--Baxter equation. Note that such monoids are YB-semitrusses $(A,+,\circ,\lambda,\id)$
as studied in \cite{CoJeVAVe21x}. For a monoid $M$ of IG-type, it is shown in \cite[Proposition~2.4]{MR2290908}
that $K[M]$ is left and right Noetherian if and only if $A$ is finitely generated as a monoid (or equivalently,
$M$ is a finitely generated monoid). Furthermore in \cite[Theorem~2.5]{MR2290908} it is shown that such monoids
are epimorphic images of monoids of I-type. We now extend this result to the structure algebra $K[M(X,r)]$ of a finite
left non-degenerate solution $(X,r)$ of the Yang--Baxter equation with abelian left derived structure monoid $A(X,r)$
(note that such a monoid is not necessarily cancellative). The proof follows the same lines as the proof of
\cite[Theorem~2.5]{MR2290908}.

\begin{prop}\label{ApplAbelian}
    Let $(X,r)$ be a finite left non-degenerate solution of the Yang--Baxter equation. Let $K$ be a field.
    If $A=A(X,r)$ is abelian then $M=M(X,r)$ is an epimorphic image of a finitely generated monoid of $I$-type
    (i.e., the structure monoid of a finite non-degenerate involutive solution of the Yang--Baxter equation).
    \begin{proof}
        Since the identity element of $A$ is the only invertible element, the generators of $A$ are unique and
        it is the set $X=\{a_1,\dotsc,a_n\}$. Let $\mc{G}=\mc{G}(X,r)=\free{\lambda_x:x\in X}$ be the permutation
        group of the solution $(X,r)$. Consider the free abelian monoid $F_m$ of rank $m=n|\mc{G}|$ with basis
        \[Y=\{y_{i,g}:1\le i\le n\text{ and }g\in\mc{G}\}.\] Since $A$ is abelian, by assumption, we obtain a monoid
        epimorphism $f\colon F_m\to A$ defined by $f(y_{i,g})=g (a_i)$. For each $\alpha\in F_m$ define the
        mapping \[\Psi_{\alpha}\colon Y \to Y\colon y_{i,g}\mapsto y_{i,\lambda_{f(\alpha)}\circ g}.\]
        Each $\Psi_{\alpha}$ is a permutation of $Y$ and thus defines an automorphism on $F_m$ (which, for simplicity,
        we also denote as $\Psi_{\alpha}$). We claim the following equalities hold for $\alpha,\beta\in F_m$:
        \begin{align}
            f(\Psi_{\alpha}(\beta)) & =\lambda_{f(\alpha)}(f(\beta)),\label{inv1}\\
            \Psi_{\alpha\Psi_{\alpha}(\beta)} & =\Psi_\alpha\circ\Psi_{\beta}.\label{inv2}
        \end{align}
        It is sufficient to check that the first equality holds for $\beta=y_{i,g}$, i.e.,
        for a basis element. This is shown as follows:
        \[f(\Psi_{\alpha}(\beta))=f(\Psi_{\alpha}(y_{i,g}))=f (y_{i,\lambda_{f(\alpha)}\circ g})
        =(\lambda_{f(\alpha)}\circ g)(a_i)=\lambda_{f(\alpha )}(f(y_{i,g}))=\lambda_{f(\alpha)}(f(\beta)).\]
        We now prove equality \eqref{inv2} by computing the image of a basis element $y_{i,g}\in Y$:
        \begin{align*}
            \Psi_{\alpha\Psi_{\alpha}(\beta)}(y_{i,g})
            & =y_{i,\lambda_{f(\alpha\Psi_{\alpha}(\beta))}\circ g}\\
            & =y_{i,\lambda_{f(\alpha)+f(\Psi_{\alpha}(\beta))}\circ g}\\
            & =y_{i,\lambda_{f(\alpha)+\lambda_{f(\alpha)}(f(\beta))}\circ g}\quad(\text{because of }\eqref{inv1})\\
            & =y_{i,\lambda_{f(\alpha)}\circ\lambda_{f(\beta)}\circ g}\\
            & =\Psi_{\alpha} (y_{i,\lambda_{f(\beta)}\circ g})\\
            & =(\Psi_{\alpha}\circ\Psi_{\beta})(y_{i,g}).
        \end{align*}
        So we get a monoid of $I$-type
        \[S=\{(\alpha,\Psi_{\alpha}):\alpha\in F_m\}\s F_m\rtimes\{\Psi_{\alpha}:\alpha\in F_m\}.\]
        Clearly \[f^e\colon S\to M\colon(\alpha,\Psi_{\alpha })\mapsto(f(\alpha),\lambda_{f(\alpha)})\]
        is a monoid epimorphism.
    \end{proof}
\end{prop}

Note that in the proof the epimorphism $f^e\colon S\to M$ is not necessarily a YB-semitruss homomorphism
for otherwise the $\sigma$-map defining $A$ would be such that $\sigma_x=\id$ for all $x\in X$. 

It is easy to construct examples of finite left non-degenerate solutions $(X,r)$ with $A(X,r)$ abelian while $M(X,r)$ is not abelian.

\begin{ex}
    Let $X=\{1,2,3\}$. Define $r\colon X\times X\to X\times X$ as $r(x,y)=(\lambda_x(y),\rho_y(x))$, where $\lambda_x=(23)$
    for each $x\in X$ and \[\rho_1(x)=1,\quad\rho_2(x)=\begin{cases}1 & \text{if }x\ne 3,\\ 2 & \text{if }x=3,\end{cases}\quad
    \rho_3(x)=\begin{cases}1 & \text{if }x\ne 2,\\ 3 & \text{if }x=2.\end{cases}\]
    Then one may check that $(X,r)$ is a left non-degenerate solution of the Yang--Baxter equation. Moreover, the monoid
    \[A(X,r)=\free{1,2,3\mid 1+1=1+2=1+3=2+1=2+3=3+1=3+2}\] is abelian, whereas
    \[M(X,r)=\free{1,2,3\mid 1\circ 1=1\circ 2=1\circ 3=2\circ 1=2\circ 2=3\circ 1=3\circ 3}\] is not abelian,
    as $2\circ 3\ne 3\circ 2$ in $M(X,r)$. 
\end{ex}

Example~\ref{ex:idem} shows that converse of the Theorem~\ref{idealchainNoetherian} does not hold, i.e., if
$K[M(X,r)]$ is left Noetherian then for each subset $Y$ of $X$ with $Y\in\mc{L}_u(i)$ for some $i$ there exists
a positive integer $d$ so that $(dw_Y,\id)\circ M_{YY}$ is cancellative. Nevertheless, we now show that for $Y=X$ the converse is true. 

\begin{prop}\label{rightnoethcancellative}
    Let $(X,r)$ be a finite left non-degenerate solution of the Yang--Baxter equation. Let $K$ be a field and $M=M(X,r)$.
    If $K[M]$ is right Noetherian or semiprime, then $M_{XX}^d$ is cancellative for some positive integer $d$.
    In particular, $M_{XX}^d$ has a group of fractions that is abelian-by-finite.
    \begin{proof}
        By Lemma~\ref{lemma8} we know that there exists a positive integer $d$ such that
        $M_{XX}^d=\bigcup_{v\in V}S_{X,v}\s\mc{M}(G,1,V,P)$, where $P$ is the $V\times 1$ matrix with all entries equal
        to $1$, the identity of $G$, and $S_{X,v}=\{(g,1,v):g\in G\} \s \mc{M}(G,1,V,P)$ are cancellative semigroups whose
        group of fractions is isomorphic with $G$. For simplicity we write elements of $S_{X,v}$ as pairs $(g,v)$ with
        $g\in G$ and $v\in V$, and thus $(g,v)\circ(g',v')=(gg',v')$.
        
        Assume $K[M]$ is semiprime, hence so is its ideal $K[M_{XX}^d]$. If $v$ and $v'$ are distinct elements of $V$,
        then choose $(g,v)\in S_{X,v}$ and $(g',v')\in S_{X,v'}$. Recall that $K[M]$ is PI by \eqref{S1}, hence so is
        each $K[S_{X,v}]$ and thus also $K[G]$. Consequently, $G$ is abelian-by-finite by \cite[Theorem~3.1.9]{JO}.
        Therefore, without loss of generality, may assume that $g$ and $g'$ commute. Hence $(g,v)\circ(g',v')=(gg',v')\in S_{X,v'}$
        and $(g',v')\circ(g,v)=(g'g,v)\in S_{X,v}$. So, replacing $g$ and $g'$ by $gg'=g'g$, we thus also may assume that $g=g'$.
        In $K[M_{XX}^d]$ it follows that $((g,v)-(g,v'))\circ M_{XX}^d=0$. Hence $((g,v)-(g,v'))\circ K[M_{XX}^d]$ is a non-zero
        nilpotent right ideal, in contradiction with the assumption. Consequently, $|V|=1$ and thus $M_{XX}^d=S_{X,v}$,
        a cancellative semigroup. 
        
        Assume now that $K[M]$ is right Noetherian. We define the following equivalence relation $\sim$ on $M$.
        For $m,m'\in M$ we say $m\sim m'$ if and only if $m=m'$ or $m=(g,v)\in S_{X,v}$ and $m'=(g',v)\in S_{X,v}$
        have the same length (as element of $M$), for some $v\in V$ and $g,g'\in G$. Clearly $\sim$ is a left
        congruence relation on $M$. To show it also is a right congruence, let $a,b\in S_{X,v}$ be elements
        of the same length (i.e., $a\sim b$) and let $m\in M$. Clearly $a\circ m,b\circ m\in M_{XX}^d$.
        Hence $a\circ m\in S_{X,v'}$ and $b\circ m\in S_{X,v''}$ for some $v',v''\in V$. Thus we need to
        prove that $v'=v''$. As the cancellative semigroup $S_{X,v}$ is a left Ore semigroup, there exist
        $a',b'\in S_{X,v}$ with $ a'\circ a=b'\circ b\in S_{X,v}$. Clearly $a'\circ a\circ m\in S_{X,v}\circ S_{X,v'}\s S_{X,v'}$
        and $b'\circ b\circ m\in S_{X,v}\circ S_{X,v''}\s S_{X,v''}$. However, as $a'\circ a\circ m=b'\circ b\circ m$,
        we get that $v'=v''$. Thus indeed $\sim$ is right congruence and thus a congruence. So, we can consider the
        semigroup $M/{\sim}$. Since, by assumption, $K[M]$ is right Noetherian also the algebra $K[M/{\sim}]$ is right Noetherian. 
        Note that $M/{\sim}$ also has an ideal chain (inherited from $M$) and remains a graded monoid (by the length function), 
        but now with lowest ideal in the ideal chain a subsemigroup of $\mc{M}(\free{g},1,V,P)$ over an infinite cyclic group
        $\free{g}$, and thus a union of cancellative semigroups that has in each cancellative component a unique element of
        a given degree (or none). Assume the lowest ideal has at least two columns, i.e., assume $|V|\ge 2$. Without loss of
        generality we may assume that $1,2\in V$. Clearly there is an infinite set $N$ of positive integers so that
        both $S_{X,1}$ and $S_{X,2}$ contain an element of degree $n\in N$. Let $s_{v,n}$ for $v\in\{1,2\}$ denote the
        unique element of degree $n\in N$ in $S_{X,v}/{\sim}$. Clearly $s_{1,n}\circ m=s_{2,n}\circ m$ for all $1\ne m\in M/{\sim}$
        and $n\in N$. Consider now the following right ideal in $K[M/{\sim}]$:
        \[\sum_{n \in N}(s_{1,n}-s_{2,n})\circ K[M/{\sim}]=\sum_{n \in N}(s_{1,n}-s_{2,n})K.\] This is an infinite dimensional
        $K$-vector space and it is thus not finitely generated as a right ideal. Hence $K[M/{\sim}]$ is not right Noetherian,
        a contradiction. So, $|V|=1$ and thus $M_{XX}^d$ is cancellative.
    \end{proof}
\end{prop}

For left derived solutions of the Yang--Baxter equation we obtain a complete characterization of when its structure
algebra is right Noetherian. It holds precisely when the converse of Theorem~\ref{idealchainNoetherian} holds.
To correctly interpret the result we give a word of caution concerning the relation between the sets $M_{YY}=M(X,r)_{YY}$
and $A_{YY}=A(X,r)_{YY}$ for a subset $Y$ of $X$. For the latter we consider $A=A(X,r)$ as the structure monoid $M(X,s)$
of the left derived solution $(X,s)$, so, in particular, the $\lambda$-map for $A=M(X,s)$ is such that $\lambda_a=\id$ for each $a\in A$.
Hence \[A_{YY}=\{a\in A:a\in A \setminus A_{|Y|+1}\text{ and }y\mid a\text{ for all }y\in Y\}.\] So, $A_{YY}=\pi (M_{Y*})$
(recall that $\pi\colon M\to A$ is the $1$-cocycle defined in Section~\ref{sec:desc}), in general this is a set of much larger cardinality
than $M_{YY}$. Also note that if $Y$ and $Z$ are different subsets of $X$ of the same cardinality then $A_{YZ}\s A_{|Y|+1}$.

\begin{cor}\label{ArightNoetherian}
    Let $(X,r)$ be a finite left non-degenerate solution of the Yang--Baxter equation. Let $K$ be a field and $A=A(X,r)$.
    Then $K[A]$ is right Noetherian if and only if for each subset $Y$ of $X$  with $(A_{YY}\cup A_{|Y|+1})/A_{|Y|+1}$
    not nil there exists
    a positive integer $d$ so that $dA_{YY}$ is cancellative.
    
    Furthermore, if $K[A]$ is right Noetherian and $M=M(X,r)$ then so is $K[M]$.
    \begin{proof}
        Let $Y$ be a subset of $X$ with $(A_{YY}\cup A_{|Y|+1})/A_{|Y|+1}$ not nil. From Lemma~\ref{lemma3} we know that the left derived
        solution $s$ of $r$ restricts to a left derived solution $s_Y$ on $Y$ and clearly the submonoid $\free{Y}$ of $A$
        generated by $Y$ satisfies $\free{Y}\cong A(Y,s_Y)$. Assume $K[A]$ is right Noetherian. Let $R=K[\free{Y}]$. Then
        consider the map $\pr\colon K[A]\to R$ defined by $\pr(\alpha+\beta)=\alpha$, where $\Supp(\alpha)\s\free{Y}$ and
        $\Supp(\beta)\s A\setminus\free{Y}$. This mapping is a right $R$-module homomorphism. Furthermore, it easily is
        verified that if $V$ is a right ideal of $R$ then $(V\cdot K[A])\cap R=V$. Hence $R\cong K[A(Y,s_Y)]$ is right Noetherian
        as well. It follows from Proposition~\ref{rightnoethcancellative} that $dA(Y,s_Y)_{YY}=dA_{YY}$ is cancellative
        for some positive integer $d$. The converse follows at once from Theorem~\ref{idealchainNoetherian} (as each
        $dA_{YY}$ is cancellative). This proves the first part of the statement.
        
        For the second part assume $K[A]$ is right Noetherian. To prove $K[M]$ is right Noetherian, because of
        Theorem~\ref{idealchainNoetherian}, it is sufficient to consider subsets $Y$ of $X$ with $M_{YY}$ a non-empty
        semigroup, $w_Y\in M_{YY}$ (as in Lemma~\ref{lemma6}) and show that $(dw_Y,\id)\circ M_{YY}$ is a cancellative
        semigroup, or equivalently high powers of elements of this semigroup commute. Clearly appropriate high powers
        of such elements are of the form $(a,\id)$ with $a\in (dw_Y+A)\cap A_{YY}$. The assumption on $Y$ also implies
        that $(A_{YY}\cup A_{|Y|+1})/A_{|Y|+1}$ is not nil. Hence $dA_{YY}$ is cancellative. Recall that $K[A]$ is a
        PI-algebra and hence so is the algebra $K[dA_{YY}]$. Because $dA_{YY}$ is a cancellative semigroup we have that
        high multiples of elements of $dA_{YY}$ commute. Consequently also high powers of the considered elements $(a,\id)$
        commute, as desired.
    \end{proof}
\end{cor}

As an application we also recover result \eqref{N3} on finite left non-degenerate idempotent solutions of the Yang--Baxter equation.

\begin{cor}\label{idempotentexample}
    Let $(X,r)$ be a finite left non-degenerate idempotent solution of the Yang--Baxter equation.
    If $K$ is a field and $M=M(X,r)$ then the following properties hold:
    \begin{enumerate}
        \item $K[M]$ is left Noetherian and $\GK K[M]=1$.
        \item $K[M]$ is right Noetherian if and only if the set $\Lambda=\{q(x)=\lambda_x^{-1}(x):x\in X\}$ is a singleton.
        The latter also is equivalent with $K[M]$ not being a central algebra, i.e, its center is strictly larger than $K$.
        \item $K[M]$ is right Noetherian if and only if there exists a positive integer $d$ such that $d$-powers of arbitrary
        elements of $M$ commute.
    \end{enumerate}
    \begin{proof}
        Let $d$ be a positive integer so that $\lambda_a^d=\id$ for each $a\in A=A(X,r)$. Clearly, since $r^2=r$, we have
        \[A=\free{X\mid x+y=y+y\text{ for all }x,y\in X}=\bigcup_{x\in X}\free{x}\] with $\free{x}\cap\free{y}=\{1\}$
        for $x\ne y$. Notice that if $a\in A$ is a word of length at least two then it is left divisible by all elements of $X$.
        Hence for a subset $Y$ of $X$ with $|Y|=i$, if $M_{YY}$ is not nil modulo $M_{i+1}$ then $Y=X$.
        One has $(x,\lambda_x)^d=(dq(x),\id) $ and thus $(x,\lambda_x)^d\circ(y,\lambda_y)^d=(y,\lambda_y)^{2d}$ for $x,y\in X$. It follows that $M_{XX}^d$
        is not a cancellative semigroup if $|\Lambda|>1$. Hence in this case, $K[M]$ is not right Noetherian by 
        Proposition~\ref{rightnoethcancellative}. On the other hand, if $\Lambda$ is a singleton then the proof shows
        that elements of the type $(a,\lambda_a)^d$ commute. Hence by Theorem~\ref{idealchainNoetherian}, $K[M]$ is right Noetherian.
        
        Clearly the description of $A$ yields that $K[A]$ has Gelfand--Kirillov dimension equal to $1$. For the last part of (2)
        we refer to \cite{MR4728712}. Note that proof also shows that $\Lambda$ is a singleton precisely when $d$-powers of arbitrary
        elements of $M$ commute. Hence the equivalence with (3) follows.
    \end{proof}
\end{cor}

Because of  Corollary~\ref{idempotentexample} it now is easy to produce finite left non-degenerate
solutions with corresponding structure algebra that is not right Noetherian.

\begin{ex}
    Consider the following left non-degenerate idempotent solution on a finite set $X$ that is not a singleton:
    $r(x,y)=(y,y)$. So, $\lambda_x=\id$ for all $x\in X$ and for each $y\in X$ we have that $\rho_y$ is a constant
    map on $X$. From Corollary~\ref{idempotentexample} we know that $K[M(X,r)]=K[A(X,r)]$ is not right Noetherian.
\end{ex}

\begin{ex}
    Let $X$ be a non-trivial finite group. Then $r(g,h)=(gh,1)$, for $g,h\in X$, defines a left non-degenerate idempotent solution
    $(X,r)$ of the Yang--Baxter equation. By Corollary~\ref{idempotentexample} the algebra $K[M(X,r)]$ is right Noetherian but
    $K[A(X,r)]$ is not right Noetherian. So, the converse of the second statement in Corollary~\ref{ArightNoetherian} does not hold. 
\end{ex}

We now show that for bijective solutions also \eqref{N3} can be obtained as an immediate application.

\begin{cor} \label{IdempotentNoetherian}
    Let $(X,r)$ be a finite left non-degenerate bijective solution of the Yang--Baxter equation.
    If $K$ is a field and $M=M(X,r)$ then $K[M]$ is left and right Noetherian.
    \begin{proof}
        Because of the bijective assumption we know, by \eqref{R1}, that $(X,r)$ also is right non-degenerate and thus
        in $A=A(X,r)$ we have $a+A=A+a$ for all $a\in A$. Hence with the notation as in Theorem~\ref{idealchainNoetherian},
        for $a,b\in dw_Y+A_{YY}$ with $ M_{YY}$ a semigroup, we have that $a+b=c+a$ for some $c\in A$. Hence also
        $(a+a)+b=(a+c)+a=(c'+a)+a$ for some $c'\in A$. So also $(a+a+a)+b=(a+c'+a)+a$ and $a+c'+a\in dw_Y+A_{YY}$.
        So, by Lemma~\ref{lemma8}, $(dw_Y,\id)\circ M_{YY}$ is a cancellative semigroup. Hence again by
        Theorem~\ref{idealchainNoetherian}, $K[M]$ is right Noetherian, as claimed. Because of the right non-degeneracy,
        the left-right dual version of this result  also yields that $K[M]$ is left Noetherian. Of course it also follows
        directly from Theorem~\ref{theorem:A}.
    \end{proof}
\end{cor}

We also easily can show that the right Noetherianity of $K[M(X,r)]$ can be reduced to that of $K[\Soc(M(X,r))]$. 

\begin{prop}\label{MtoSoc}
    Let $(X,r)$ be a finite left non-degenerate solution of the Yang--Baxter equation.
    If $K$ is a field, $M=M(X,r)$ and $S=\Soc(M)$ then $K[M]$ is right Noetherian if and only if $K[S]$ is right Noetherian.
    \begin{proof}
        As Proposition~\ref{fromMtoS} assures that $K[M]$ is a finite right module over $K[S]$, it remains to show that $K[S]$
        is right Noetherian if $K[M]$ is so. Therefore, assume that $K[M]$ is right Noetherian. Since $(s\circ M)\cap S=s\circ S$
        for each $s\in S$, it easily follows that $(V\circ K[M])\cap K[S]=V$ for any right ideal $V$ of $K[S]$. This clearly
        implies that $K[S]$ is right Noetherian.
    \end{proof}
\end{prop}

The ideal chain obtained in Theorem~\ref{idealchainNoetherian} intersected with $S=\Soc(M)$ yields an ideal chain for
$S$ with the same properties. Proceeding with a proof as in Corollary~\ref{ArightNoetherian} and making use of
Lemma~\ref{lemma3} we obtain at once a criterion for $K[S]$ to be right Noetherian. Hence Proposition~\ref{MtoSoc}
leads immediately to the following characterization when $K[M]$ is right Noetherian. But before formulating it, we shall
introduce the following notation. For non-empty subsets $Y,Z$ of $X$ of the same cardinality let \[S_{YZ}=M_{YZ}\cap S.\]
Note that if $Y\ne Z$ then $S_{YZ}=\vn$.

\begin{thm}\label{SrightNoetherian}
    Let $(X,r)$ be a finite left non-degenerate solution of the Yang--Baxter equation. If $K$ is a field, $M=M(X,r)$
    and $S=\Soc(M)$ then $K[M]$ is right Noetherian if and only if there exists a positive integer $d$ so that
    $S_{YY}^d$ is cancellative for each subset $Y$ of $X$ with $(S_{YY}\cup M_{|Y|+1})/M_{|Y|+1}$ not nil.
\end{thm} 

\section{Prime structure algebras and cancellative structure monoids}\label{sec:prime}

Let $(X,r)$ be a finite left non-degenerate solution of the Yang--Baxter equation.
Recall \cite{ESS99} the definition of its associated structure group:
\begin{align*}
    G(X,r) & =\gr(x_1,\dotsc,x_n\mid x_i \circ x_j=x_k \circ x_l\text{ if }r(x_i,x_j)=(x_k,x_l))\\
    & =\gr(x_1,\dotsc,x_n\mid x_i \circ x_j=\lambda_{x_i}(x_j)\circ\rho_{x_j}(x_i)).
\end{align*}
Here we denote by $\gr(x_1,\dotsc,x_n\mid u_1=v_1,\dotsc,u_m=v_m)$ the group defined by free generators
$x_1,\dotsc,x_n$ with respect to the relations $u_i=v_i$ for $1\le i\le m$, where $u_i,v_i$ are words in
the free monoid $\free{x_1,\dotsc,x_n}$. There is an obvious natural monoid homomorphism $M(X,r)\to G(X,r)$.
In general this is not injective (the images of elements of $X$ could become equal in $G(X,r)$, see, e.g.,
\cite{MR1809284}). Of course, if $M(X,r)$ has a group of fractions then this group is isomorphic with $G(X,r)$.

\begin{prop}\label{prime}
    Let $(X,r)$ be a finite left non-degenerate solution of the Yang--Baxter equation. Let $K$ be a field and $M=M(X,r)$.
    If $K[M]$ is a prime algebra then $M$ is a cancellative monoid with group of fractions isomorphic to the structure group $G=G(X,r)$. 
    
    Furthermore, if $S=\Soc(M)$ then $S\circ S^{-1}$ is a free abelian group of finite rank contained in the finite conjugacy
    center of $G$, $S$ is a left and right Ore subset of $M$ and $G=M\circ S^{-1}=M\circ\Z(M)^{-1}$. In particular,
    $K[G]$ is a central localization of $K[M]$ and thus $K[G]$ is a prime algebra.
    \begin{proof}
        From Proposition~\ref{rightnoethcancellative} we know that the ideal $M_{XX}^d$ is cancellative for some $d\ge 1$
        and it has an abelian-by-finite group of quotients, say $U$. Since the algebra $K[M]$ is prime and PI, it follows
        by \cite[Corollary~13.6.6]{MR} that $K[M]$ is a (two-sided) Goldie ring with simple Artinian classical ring of
        quotients $Q=\Q_{cl}(K[M])$. By Posner's theorem (see, e.g., \cite[Theorem~13.6.5]{MR}) and the fact that central elements
        of the ideal $K[M_{XX}^d]$ are central in the prime algebra $K[M]$, we get that $Q=K[M]\circ Z^{-1}$, a central
        localization of $K[M]$ with respect to the multiplicative set
        \[Z=(\Z(K[M])\cap K[M_{XX}^d])\setminus\{0\}=\Z(K[M_{XX}^d])\setminus\{0\}.\] Therefore, $K[M_{XX}^d]\circ Z^{-1}$
        is a non-zero ideal of $Q$ and thus $K[M_{XX}^d]\circ Z^{-1}=Q$, as $Q$ is a simple ring.
        Because $U$ is the group of quotients of $M_{XX}^d$, we get that $\Q_{cl}(K[U])=K[M_{XX}^d]\circ Z^{-1}=Q$.
        
        To prove that the monoid $M$ is cancellative, assume $a\circ c=b\circ c$ for some $a,b,c\in M$. Then taking
        $s\in M_{XX}^d$ we get $(s\circ a)\circ (c\circ s)=(s\circ b)\circ (c\circ s)$. Since $s\circ a,s\circ b,c\circ s\in M_{XX}^d$
        we get $s\circ a=s\circ b$ and thus $a=s^{-1}\circ (s\circ a)=s^{-1}\circ (s\circ b)=b$ in $\Q_{cl}(K[U])=Q$.
        As $M$ embeds in $Q$, we conclude that $a=b$ in $M$. Similarly, one shows that $c\circ a=c\circ b$ implies $a=b$.
        Therefore, the monoid $M$ is cancellative. This finishes the proof of the first part of the statement.
        
        That $G=M\circ\Z(M)^{-1}$ follows at once from \cite[Lemma~4.1.9]{JO}. For completeness' sake we include a proof.
        Let $g\in G$. Since $g\in\Q_{cl}(K[M])$, there exists a central element $\alpha\in K[M]$ so that $g\circ \alpha\in K[M]$.
        Hence $g\circ\Supp(\alpha)\s M$. Now, for any $h\in G$, we have $h\circ\alpha\circ h^{-1}=\alpha$ and thus
        $h\circ\Supp(\alpha)\circ h^{-1}=\Supp(\alpha)\s M$. Therefore, $h\circ x\circ h^{-1}\in M$ and $g\circ h\circ x\circ h^{-1}\in M$
        for any $x\in\Supp(\alpha)$. Since $G$ is abelian-by-finite, it has an abelian normal subgroup $V$ of finite index, say $t$.
        In particular $x^t\in V\cap M$ and $x^t$ has only finitely many conjugates in $G$, say $g_i\circ x^t\circ g_i^{-1}$ with
        $1\le i \le m$, and all these conjugates commute. Thus we get that
        \[z=(g_1\circ x\circ g_1^{-1})^t\circ\dotsb\circ(g_m\circ x\circ g_m^{-1})^t\in\Z(M).\]
        Then  $g\circ z=g\circ (g_1\circ x\circ g_1^{-1})^t\circ\dotsb\circ(g_m\circ x\circ g_m^{-1})^t\in M$ and
        consequently $g\in M\circ\Z(M)^{-1}$, as claimed.
        
        To prove the second part recall that from Proposition~\ref{fromMtoS}(2) we know that $S=\Soc(M)$ is finitely generated.
        Since certain high powers of elements of $G$ and thus also appropriate high powers of $M$ commute and belong to $S$,
        it follows that $S$ is a left and right Ore subset of $M$ and $G=M\circ S^{-1}$. Put $H=S\circ S^{-1}$, a subgroup of $G$.
        Now, for $(s,\id)\in S$ and $(a,\lambda_a)\in M$ we have
        \[(s,\id)\circ(a,\lambda_a)=(a,\lambda_a)\circ(\lambda_a^{-1}(\sigma_a(s)),\id)\]
        and thus in $G$: \[(a,\lambda_a)^{-1}\circ(s,\id)\circ(a,\lambda_a)=(\lambda_a^{-1}(\sigma_a(s)),\id).\]
        Because $\lambda_a^{-1} \sigma_a$ is length preserving and $G=M\circ\Z(M)^{-1}$, it follows that each element
        of $S$ has only finitely many conjugates in $G$ and, furthermore, $H$ is a normal subgroup of $G$.
        In particular, $H$ is contained in the finite conjugacy center of $G$. In particular, $H$ is a finitely
        generated finite conjugacy group. It is well known that the set $T(H)$ consisting of the periodic elements of 
        $H$ is a then a finite subgroup of $H$ and $H/T(H)$ is a finitely generated free abelian group. As $T(H)$
        is a characteristic subgroup of $H$ and $H$ is a normal subgroup of $G$, we get that $T(H)$ is a normal
        subgroup of $G$. Therefore, \[\gamma=\sum_{h\in T(H)}h\] is a central element of $K[G]$ with the property
        that $\gamma^2=|T(H)|\gamma$. Hence $\gamma(\gamma-|T(H)|1)=0$. Since $\gamma\ne 0$ and $K[G]$ is prime,
        it follows that $\gamma=|T(H)|1$ and thus $T(H)=\{1\}$. Consequently, $H$ is finitely generated torsion-free
        abelian group, thus a finitely generated  free abelian group.
    \end{proof}
\end{prop}

Since the semigroup $\free{\sigma_x:x\in X}$ is finite there exists a positive integer, say $e_{\sigma}$
so that each $\sigma_x^{e_\sigma}$ is an idempotent mapping, i.e.,
\begin{equation}\label{idempotentexponent}
    \sigma_x^{2e_{\sigma}}=\sigma_x^{e_{\sigma}}.
\end{equation}
Clearly we may replace $e_{\sigma}$ by any multiple. Let
\begin{equation}\label{defv}
    v=kee_{\sigma},
\end{equation}
where $k$ and $e$ are the positive integers introduced in \eqref{ekdef} and \eqref{expdef}. From now on we replace
the set $W$ defined in \eqref{defW} by $W_v$, i.e., \[W=W_v=\{(vx,\id):x\in X\text{ and }\lambda_{vx}=\id\}.\]
Again let us write $X=\{x_1,\dotsc,x_n\}$. For $\kappa\in\Sym(n)$ put
\begin{equation}
    z_\kappa=vx_{\kappa(1)}+\dotsb+vx_{\kappa(n)}\in A.\label{defz}
\end{equation}
In \cite[Equation~(5.6)]{CoJeVAVe21x} it has been shown that if $x,y\in X$ then 
\begin{equation}
    vx+vy+z_\kappa=vy+vx+z_\kappa.\label{zcomm}
\end{equation}
It is worth mentioning that in the proof of \eqref{zcomm} the fact that $\sigma_{vx},\sigma_{vy},\sigma_{z_\kappa}$
are idempotent maps is crucial. Moreover, this property is essential in the proof result \eqref{N1}.

\begin{prop}\label{primeAcancellative}
    Assume $(X,r)$ is a finite left non-degenerate solution of the Yang--Baxter equation.
    Let $K$ be a field and let $W=W_v$. If the monoid $M=M(X,r)$ is cancellative then:
    \begin{enumerate}
        \item $A=A(X,r)$ is a left cancellative monoid and $S=\Soc(M)$ is a cancellative monoid.
        \item $\gr(S)=S\circ S^{-1}$ and $\gr(W)=\free{W}\circ\free{W}^{-1}$ are normal subgroups of $G=G(X,r)$ contained in
        the finite conjugacy center of $G$. In particular, the torsion elements of these groups are finite normal subgroups of $G$.
        Furthermore, if $K[M]$ is a prime algebra then $\gr(W)$ and $\gr(S)$ are free abelian groups.
        \item $\free{W}$ is an abelian monoid.
        \item $M=\bigcup_{f\in F} f\circ\free{W}$ and $G=\bigcup_{f\in F} f\circ\gr(W)$.
        \item $K[M]$ is right and left Noetherian PI-algebra.
        \item if the diagonal map $q\colon X\to X$ is bijective then $z_{\id}=z_\kappa$ for any $\kappa\in\Sym(n)$.
        In particular, $\lambda_a(z_{\id})=z_{\id}$ for any $a\in A$.
        \item if the diagonal map $q$ is bijective then $A$ is cancellative.
    \end{enumerate}
    \begin{proof}
        (1) To prove that $A$ is left cancellative assume $a,b,c\in A$ are such that $a+b=a+c$. Then, in $M$,
        \[(a,\lambda_a)\circ(\lambda_a^{-1}(b),\lambda_{\lambda_a^{-1}(b)})=(a,\lambda_a)\circ(\lambda_a^{-1}(c),\lambda_{\lambda_a^{-1}(c)}).\] 
        Because, by assumption, $M$ is (left) cancellative, we get that
        $(\lambda_a^{-1}(b),\lambda_{\lambda_a^{-1}(b)})=(\lambda_a^{-1}(c),\lambda_{\lambda_a^{-1}(c)})$.
        The bijectivity of $\lambda_a$ therefore yields $b=c$, as desired. Since $S$ is isomorphic with a submonoid of $A$,
        it follows that it also is left cancellative. To prove it is right cancellative, assume that $a,b,c\in S$ satisfy
        $b+a=c+a$. Then, in $M$, we have $(b,\id)\circ(a,\id)=(c,\id)\circ(a,\id) $. Hence the (right) cancellativity of
        $M$ yields $b=c$, again as desired.
        
        (2) That the subgroup $\gr(S)=S\circ S^{-1}$ is normal in $G$, and contained in
        the finite conjugacy center of $G$, has been shown in the proof of Proposition~\ref{prime}, and it also shows
        that $\gr(W)=\free{W}\circ\free{W}^{-1}$ is a normal subgroup of $G$. Furthermore, if the algebra $K[M]$ is prime,
        then Proposition~\ref{prime} guarantees that $K[G]$ is a central localization of $K[M]$. In particular, $K[G]$ is prime.
        Hence \cite[Theorem~4.2.10]{MR798076} implies that the finite conjugacy center of $G$ is torsion-free and abelian.
        So, $\gr(W)$ and $\gr(S)$ are free abelian groups.
        
        (3) Suppose $(vx,\id) ,(vy,\id)\in W$. Then, by \eqref{zcomm}, $vx+vy+z_{\id}=vy+vx+z_{\id}$. Hence,
        \[(vx,\id)\circ(vy,\id)\circ(z_{\id},\lambda_{z_{\id}})=(vy,\id)\circ(vx,\id)\circ(z_{\id},\lambda_{z_{\id}}).\]
        Because $M$ is right cancellative it follows that $(vx,\id)\circ(vy,\id)=(vy,\id)\circ(vx,\id) $,
        and thus indeed $\free{W} $ is abelian.
        
        (4) From Proposition~\ref{fromMtoS}(1) we know that $M=\bigcup_{f\in F}f\circ\free{W}$. Let $f\in F$. Because $F$
        is finite, some distinct powers of $f$ lie in the same coset. Therefore, there exist positive integers $n_1,n_2$
        and $f_1\in F$ so that $f^{n_1},f^{n_1+n_2}\in f_1\circ\free{W}$. Write $f^{n_1}=f_1\circ w_1$ and $f^{n_1+n_2}=f_1\circ w_2$
        for some $w_1,w_2\in\free{W}$. By part (2), $\free{W}$ is abelian and thus $f^{n_1}\circ w_2=f^{n_1+n_2}\circ w_1$. Hence 
        by the cancellativity of $M$, $f^{n_2}\circ w_1=w_2$ and thus $f^{n_2}\in\gr(W)$. So $\gr(W)$ is a normal and abelian
        subgroup of finite index in $G$. It follows that $G=\bigcup_{f\in F}f\circ\gr(W)$.
        
        (5) Since $W$ is finite and $\free{W}$ is abelian by part (2), the commutative algebra $K[\free{W}]$ is Noetherian.
        Hence by part (3), $K[M]$ is a right Noetherian PI-algebra. From Theorem~\ref{theorem:A} we know that $K[M]$ is left Noetherian.
        
        (6) Let $\kappa\in\Sym(n)$. Because $q$ is bijective, we get $W=W_v=\{(vx,\id):x\in X\}$, which yields
        $(z_\kappa,\id)\in\free{W}\s\Soc(M)$. Since, by \eqref{zcomm}, $z_\kappa+z_{\id}=z_{\id}+z_{\id}$ then
        $(z_\kappa,\id)\circ(z_{\id},\id)=(z_{\id},\id)\circ(z_{\id},\id)$. Because $\Soc(M)$ is cancellative by part (1),
        we conclude that $z_\kappa=z_{\id}$.
        
        (7) From part (1) we already know that $A$ is left cancellative. To prove it also is right cancellative,
        suppose $a+c=b+c$ for some $a,b,c\in A$. Applying \eqref{sumsigma} several times, we may assume, without
        loss of generality, that $c=z_\kappa$ for some $\kappa\in\Sym(n)$ (here, if needed, we have to replace
        $v$ by some multiple). By part (6) $\lambda_a^{-1}(c)=z_{\id}$ for all $a\in A$. Hence in $M$, we have
        \[(a+c,\lambda_{a+c})=(a,\lambda_a)\circ(\lambda_a^{-1}(c),\lambda_{\lambda_a^{-1}(c)})=(a,\lambda_a)\circ(z_{\id},\id)\]
        and \[(b+c,\lambda_{b+c})=(b,\lambda_b)\circ(\lambda_b^{-1}(c),\lambda_{\lambda_b^{-1}(c)})=(b,\lambda_b)\circ(z_{\id},\id).\]
        Because $M$ is cancellative we obtain $a=b$, as desired.
    \end{proof}
\end{prop}

We are now in a position to extend result \eqref{S1} to left non-degenerate solutions that are not
necessarily bijective provided the diagonal map $q\colon X\to X$ is bijective.

\begin{thm}\label{bijectiveqprime}
    Assume $(X,r)$ is a finite left non-degenerate solution of the Yang--Baxter equation. Let $K$ be a field
    and let $M=M(X,r)$. If the diagonal map $q\colon X\to X$ is bijective then following conditions are equivalent:
    \begin{enumerate}
        \item $(X,r)$ is an involutive solution.
        \item $M$ is a cancellative monoid.
        \item $K[M]$ is a prime algebra.
        \item $K[M]$ is a domain.
    \end{enumerate}
    \begin{proof} 
        The implication $(1)\Longrightarrow(4)$ is well-known (see \cite[Corollary~1.5]{GVdB}) and the implication
        $(4)\Longrightarrow(3)$ is trivial. Next, the implication $(3)\Longrightarrow(2)$ follows from Proposition~\ref{prime}.
        Finally, suppose than $M$ is a cancellative monoid. Then, Proposition~\ref{primeAcancellative} yields that
        $A=A(X,r)$ is cancellative and thus the left derived solution $(X,s)$ of $(X,r)$ is bijective. Consequently,
        the solution $(X,r)$ is bijective as well (see, e.g., \cite[Proposition~2.4]{CoJeVAVe21x}). Now, \cite[Theorem~4.5]{JKVA2018}
        assures that the solution $(X,r)$ is involutive. Hence $(2)\Longrightarrow(1)$, and the proof is finished.
    \end{proof}
\end{thm}

It remains an open problem to deal with the case that $q$ is not bijective. Nevertheless, in the remainder of this section,
we show that the equivalence of (4) and (1) remains valid even if $q$ is not bijective. In order to do so, we will construct
some finite subgroups in the structure group. One of the issues is that the elements $z_\kappa$ defined in \eqref{defz} do
not necessarily belong to $\Soc(M)$ (they do if $q$ is bijective). Of course $(z_\kappa,\lambda_{z_\kappa})^e\in\Soc(M)$;
recall that $e$ is the exponent of $\mc{G}(X,r)$. These elements will play a role in the construction of finite subgroups.
For simplicity, we also will denote the diagonal map of the solution $(M,r_M)$ as $q$. Clearly this is an extension
of the diagonal map of $(X,r)$. For $\kappa\in\Sym(n)$ put \[O_{\lambda}(eq(z_\kappa))=\{\lambda_a(eq(z_\kappa)):a\in A\},\]
i.e., the orbit of $eq(z_\kappa)$ under the action of the permutation group $\mc{G}(X,r)=\{\lambda_a:a\in A\}$.
Put \[T=\free{z_\kappa:\kappa\in\Sym(n)},\] a submonoid of $A$.

\begin{lem}\label{orbitprop}
    Let $(X,r)$ be a finite left non-degenerate solution of the Yang--Baxter equation. If the monoid $M=M(X,r)$ is cancellative then:
    \begin{enumerate}
        \item $eq(z_\kappa)=eq(z_{\id})$ for all $\kappa\in\Sym(n)$.
        \item for any $z\in\free{T}$ with $|z|=|eq(z_{\id})|$ and $\lambda_z=\id$ we have $z=eq(z_{\id})$.
        \item for any $z\in O_{\lambda} (eq(z_{\id}))$ and any $t\in\free{T} $ with $|t|=|z|$ we have $\sigma_z(t)=z$.
        \item $(M,r_M)$ restricts to an idempotent left non-degenerate solution on $O_{\lambda}(eq(z_{\id}))$,
        denoted by $r_{\lambda}$, and $r_{\lambda}(z,t)=(\lambda_z(t), q(\lambda_z(t)))$.
    \end{enumerate}
    \begin{proof}
        (1) Clearly
        \[(z_\kappa,\lambda_{z_\kappa})^e=(z_\kappa+\lambda_{z_\kappa}(z_\kappa)+\dotsb+\lambda_{z_\kappa}^{e-1}(z_\kappa),\id)
        =(e\lambda_{z_\kappa}^{-1}(z_\kappa),\id)=(eq(z_\kappa),\id).\]
        Further, from \eqref{zcomm} we know that
        $eq(z_\kappa)+eq(z_{\id})=eq(z_{\id})+eq(z_{\id})$. Because $\lambda_{eq(z_\kappa)}=\lambda_{eq(z_{\id})}=\id$,
        the cancellativity of $M$ implies that $eq(z_\kappa)=eq(z_{\id})$.
        
        (2) We have $z+eq(z_{\id})=eq(z_{\id})+eq(z_{\id})$. As $\lambda_{z}=\lambda_{eq(z_{\id})}=\id$,
        it follows that \[(z,\id)\circ(eq(z_{\id}),\id)=(eq(z_{\id}))\circ(eq(z_{\id}),\id).\]
        The cancellativity of $M$ thus yields $z=eq(z_{\id})$.
        
        (3) Let $t\in\free{T}$ with $|t|=|eq(z_{\id})|$. Again from \eqref{zcomm} we have
        \[eq(z_{\id})+eq(z_{\id})=t+eq(z_{\id})=eq(z_{\id})+\sigma_{eq(z_{\id})}(t).\]
        As $\lambda_{eq(z_{\id})}=\id$, it follows that $\lambda_{\sigma_{eq(z_{\id})}(t})=\id$ and 
        \[(eq(z_{\id}),\id)\circ(eq(z_{\id}),\id)=(eq(z_{\id}),\id)\circ(\sigma_{eq(z_{\id})}(t),\id).\]
        The cancellativity of $M$ yields $eq(z_{\id})=\sigma_{eq(z_{\id})}(t)$. For any $a\in A$ we know that 
        $\sigma_{\lambda_a (eq(z_{\id}))}(\lambda_a(t))=\lambda_a(\sigma_{eq(z_{\id})}(t))$. Hence
        $\sigma_{\lambda_a (eq(z_{\id}))}(\lambda_a(t))=\lambda_a (eq(z_{\id}))$. As the subsets of $\free{T}$
        of a particular length are $\lambda$-invariant we obtain that $\sigma_{\lambda_a(eq(z_{\id}))}(t)=\lambda_a(eq(z_{\id}))$
        for any $t\in\free{T}$ with $|t|=|eq(z_{\id})|$. In particular, for any $z\in O_{\lambda}(eq(z_{\id}))$
        and $t\in\free{T}$ with $|z|=|t|=|eq(z_{\id})|$, we get that $\sigma_z(t)=z$.
        
        (4) Let $t,z\in O_{\lambda}(eq(z_{\id}))$. Denote $z'=\lambda_z(t)\in O_{\lambda}(eq(z_{\id}))$.
        Because of (3), $\sigma_{z'}(z)=z'$ and thus \[r^2_M(z,t)=r_M(z',\rho_t(z))=(\lambda_{z'}(\rho_t(z)),w)=(\sigma_{z'}(z),w)=(z',w)\]
        for some $w\in M$. Consequently, $\lambda_{z'}(\rho_t(z))=z'$, which yields $\rho_t(z)=\lambda_{z'}^{-1}(z')\in O_{\lambda}(ez_{\id})$.
        
        Thus \[r_{\lambda}(z,t)=(\lambda_z(t),q(\lambda_z(t)))\in O_{\lambda}(eq(z_{\id}))\times O_{\lambda}(eq(z_{\id})).\]
        So, indeed, $r_M$ restrict to a solution on $O_{\lambda}(eq(z_{\id}))$. Also note that this restricted solution $r_{\lambda}$
        is left non-degenerate, because $O_{\lambda}(eq(z_{\id}))$ is finite and each $\lambda_z$ is injective, and thus bijective on
        $O_{\lambda}(eq(z_{\id}))$.
        
        It remains to show that the restricted solution is idempotent. For this assume $z,t\in O_{\lambda}(eq(z_{\id}))$. Then, 
        $r_{\lambda}(z,t)=(\lambda_z(t),q(\lambda_z(t)))$ and taking into account that $\lambda_a(q(a))=a$ for each $a\in A$, we obtain
        \begin{align*}
            r_{\lambda}^2(z,t) &=r_{\lambda}(\lambda_z(t),q(\lambda_z(t)))\\
            & =(\lambda_{\lambda_z(t)}(q(\lambda_z(t))),q(\lambda_{\lambda_z(t)}(q(\lambda_z(t)))))\\
            & =(\lambda_z(t),q(\lambda_z(t)))=r_{\lambda}(z,t).\qedhere
        \end{align*}
    \end{proof}
\end{lem}

\begin{prop}\label{ogroup}
    Let $(X,r)$ be a finite left non-degenerate solution of the Yang--Baxter equation. If the monoid $M=M(X,r)$ is cancellative
    then $(O_{\lambda}(eq(z_{\id})),\bullet)$ is a group for the operation $\bullet$ defined as $z\bullet t=\lambda_z(t)$,
    with identity element $eq(z_{\id})$.
    \begin{proof}
        First we prove the operation $\bullet$ is associative. So, let $z,t,t'\in O_{\lambda}(eq(z_{\id}))$. Then
        \[z\bullet(t\bullet t')=z\bullet\lambda_t(t')=\lambda_z(\lambda_t(t'))=\lambda_{\lambda_z(t)}(\lambda_{\rho_t(z)}(t'))\]
        and \[(z\bullet t)\bullet t'=\lambda_z(t)\bullet t'=\lambda_{\lambda_z(t)}(t').\] From Lemma~\ref{orbitprop}(4) we know
        that $\rho_t(z)=q(\lambda_z(t))$. Hence for associativity of $\bullet$ we need to prove that 
        $\lambda_{\rho_t(z)}=\lambda_{q(\lambda_z(t))}=\id$ and this follows from the following:
        \begin{align*}
            (\lambda_z(t),\lambda_{\lambda_z(t)})^e & =(\lambda_z(t)+\dotsb+q(\lambda_z(t)),\id)\\
            & =((e-1)eq(z_{\id})+q(\lambda_z(t)),\id)\\
            & =(eq(z_{\id}),\id)^{e-1}\circ (q(\lambda_z(t)),\lambda_{q(\lambda_z(t))}).
        \end{align*}
        Because $\lambda_{q(\lambda_z(t))}=\id$ and $|q(\lambda_z(t))|=|eq(z_{\id})|$ we obtain from Lemma~\ref{orbitprop}(2)
        that $q(\lambda_z(t))=eq(z_{\id})$. In particular, for $z=eq(z_{\id})$ we get that $q(t)=eq(z_{\id})$, i.e.,
        $t=\lambda_t(eq(z_{\id}))$. 
        
        Clearly, $\lambda_{eq(z_{\id})}=\id$ implies that $eq(z_{\id})\bullet t=t$ and thus $eq(z_{\id})$ is a left identity.
        That $eq(z_{\id})$ also is a right identity  follows from $t\bullet eq(z_{\id})=\lambda_t(eq(z_{\id}))=t$.
        
        Finally, $z\bullet\lambda_z^{-1}(eq(z_{\id}))=eq(z_{\id})$ and thus every element $z\in O_{\lambda}(eq(z_{\id}))$
        is right invertible. Hence indeed $O_{\lambda}(qe(z_{\id})$ is a finite group with $eq(z_{\id})$ as the identity element.
    \end{proof}
 \end{prop}

\begin{prop}\label{finitesubgroup}
    Let $(X,r)$ be a finite left non-degenerate solution of the Yang--Baxter equation. If the monoid $M(X,r)$ is cancellative then
    \[\Omega_\lambda=\{(z,\lambda_z)\circ(eq(z_{\id}),\id)^{-1}:z\in O_{\lambda}(eq(z_{\id}))\}\]
    is a finite subgroup of $(G(X,r),\circ)$.
    \begin{proof}
        First note that for $z\in O_{\lambda}(eq(z_{\id}))$ we have that
        $(z,\lambda_z)\circ (eq(z_{\id}),\id)^{-1}=(eq(z_{\id}),\id)^{-1}\circ (z,\lambda_z)$.
        Indeed, we need to show that $(z,\lambda_z)\circ(eq(z_{\id}),\id)=(eq(z_{\id}),\id)\circ(z,\lambda_z)$.
        Equivalently, we have to show that $z+\lambda_{z}(eq(z_{\id}))=eq(z_{\id})+z$ in $A$. By Proposition~\ref{ogroup},
        $\lambda_z (eq(z_{\id}))=z\bullet(eq(z_{\id}))=z$. So we need to prove that $z+z=eq(z_{\id})+z$.
        But this follows from \eqref{zcomm}.
        
        Now we prove that the set $\Omega_\lambda$ is multiplicatively closed (and thus a finite group, as it is cancellative).
        Indeed, for $z,w\in O_{\lambda}(eq(z_{\id}))$, we have
        $z+\lambda_z(w)=eq(z_{\id})+\lambda_z(w)$ by \eqref{zcomm}. Hence
        $(z,\lambda_z)\circ (w,\lambda_w)=(eq(z_{\id}),\id)\circ(\lambda_{z}(w),\lambda_{\lambda_{z}(w)})$ and thus
        \begin{align*}
            (z,\lambda_z)\circ(eq(z_{\id}),\id)^{-1}\circ(w,\lambda_w)\circ(eq(z_{\id}),\id)^{-1}
            & =(eq(z_{\id}),\id)^{-2}\circ(eq(z_{\id}),\id)\circ(\lambda_{z}(w),\lambda_{\lambda_{z}(w)})\\
            & =(eq(z_{\id}),\id)^{-1}\circ(\lambda_{z}(w),\lambda_{\lambda_{z}(w)})\in \Omega_\lambda\\
            & =(\lambda_{z}(w),\lambda_{\lambda_{z}(w)})\circ(eq(z_{\id}),\id)^{-1}\in\Omega_\lambda.\qedhere
        \end{align*}
    \end{proof}
\end{prop}

\begin{thm}\label{primefinitegroup}
    Let $(X,r)$ be a finite left non-degenerate solution of the Yang--Baxter equation.
    If $K$ is a field and $M=M(X,r)$ then the following conditions are equivalent:
    \begin{enumerate}
        \item $(X,r)$ is an involutive solution.
        \item $M$ is a cancellative monoid and $\Omega_\lambda$ is a trivial group.
        \item $K[M]$ is a prime algebra and $\Omega_\lambda$ is a trivial group.
        \item $K[M]$ is a domain.
    \end{enumerate}
    \begin{proof}
        If $K[M]$ is a domain then, of course, $K[M]$ is a prime algebra. From Proposition~\ref{prime} we then also get
        that $M$ is cancellative and $G=G(X,r)$ is its group of fractions, so also $K[G]$ is a domain. Consequently,
        $G$ has no non-trivial elements of finite order. Hence by Proposition~\ref{finitesubgroup}, $\Omega_\lambda$
        is trivial, i.e., $O_{\lambda}(eq(z_{\id}))=\{eq(z_{\id})\}$, that is $eq(z_{\id})$ is $\lambda$-invariant.
        It then follows easily that the left derived monoid $A=A(X,r)$ is right cancellative and thus cancellative.
        Hence the remaining part of the proof is as in the proof of Theorem~\ref{bijectiveqprime}.
    \end{proof}
\end{thm}

Because of Theorem~\ref{bijectiveqprime} the condition $|\Omega_\lambda|=1$ in the statement is redundant provided the
diagonal map is bijective. It remains an open problem whether $\Omega_\lambda$ is a singleton when $K[M]$ is prime.

\section{Semiprime Noetherian structure algebras}\label{sec:semiprime}

Let $(X,r)$ be a finite non-degenerate  bijective solution of the Yang--Baxter equation. For the left derived solution
$(X,s)$ the diagonal map $q$ is bijective. Hence from Theorem~\ref{primefinitegroup} we obtain that $K[A]$ is prime
if and only if $A=A(X,r)$ is a free abelian monoid. In this section we characterize when $K[A]$ is semiprime in terms
of a decomposition of $A$ as a finite semilattice of certain cancellative subsemigroups. 

Recall from Corollary~\ref{IdempotentNoetherian} that this algebra is left and right Noetherian.
Also recall that a semilattice is an abelian semigroup consisting of idempotents.

\begin{thm}\label{semiprimeA}
    Let $(X,r)$ be a finite non-degenerate bijective solution of the Yang--Baxter equation.
    If $K$ is a field and $A=A(X,r)$ then the following properties are equivalent:
    \begin{enumerate}
        \item the algebra $K[A]$ is semiprime,
        \item the monoid $A$ is a disjoint union $A=\bigcup_{e\in\Gamma}A_e$ of cancellative semigroups $A_e$
        indexed by a finite semilattice $\Gamma$ and $A_eA_f\s A_{ef}$ for all $e,f\in\Gamma$ (i.e., $A$ is a
        finite semilattice $\Gamma$ of cancellative semigroups) each $K[A_e]$ is semiprime.
    \end{enumerate}
    Moreover, in case the above equivalent conditions hold, $\Gamma $ is the set of central idempotents of the classical ring
    of quotients of $K[A]$. Equivalently, $A$ has a finite ideal chain with Rees factors cancellative semigroups that yield
    semiprime semigroup algebras. The latter condition holds in case $K$ has zero characteristic.
    \begin{proof}
        Assume first that the algebra $K[A]$ is semiprime. Since also  $K[A]$ is  Noetherian, it has a classical ring of
        fractions, denoted $Q=\Q_{cl}(K[A])$. Because of the bijective assumption, the left derived solution $(X,s)$ also is
        bijective (see, e.g., \cite{JKVA2018}) and thus the  elements of $A$ are normal, that is $aA=Aa$ for each $a\in A$
        (to avoid confusion with the addition in $Q$ we denote the monoid $A$ in multiplicative notation). Hence,
        $Qx$ is a non-zero ideal in the semisimple ring $Q$, for each  $x\in X$. Therefore, there exists a central
        idempotent $e_x\in Q$ such that $Qx=Qe_x$. Let $\Gamma$ be the subsemigroup of the multiplicative semigroup
        of the semisimple ring $Q$ generated by the set $\{e_x:x\in X\}$. Clearly $\Gamma$ is a finite semilattice.
        For $e\in \Gamma$ define $A_e=\{a\in A:Qa=Qe\}$. It is easy to see that each $A_e$ is a subsemigroup of $A$.
        Moreover, if $a\in A_e$ then $(Qe)a=(Qa)e=(Qe)e=Qe$, so $A_e$ consists of invertible elements of the ring
        $Qe$ and thus $A_e$ is a cancellative semigroup. Also, $A_eA_f\s A_{ef}$ for all $e,f\in\Gamma$. Furthermore,
        if $a\in A$ then $Qa$ is an ideal of the semisimple ring $Q$, and thus $Qa=Qe$ for some $e\in\Gamma$, so $a\in A_e$.
        Therefore, $A=\bigcup_{e\in \Gamma}A_e$ and it also is  clear that this is a disjoint union. Finally, if $e\in\Gamma$
        then $QA_e=Qe$. Hence the algebra $K[A_e]$ is semiprime as each nilpotent ideal of $K[A_e]$ extends to a nilpotent
        ideal of the semisimple ring~$Qe$.
        
        The reverse implication is easy (see, e.g., \cite[Corollary~3.2.9]{JO}). Also recall that $K[A_e]$
        is semiprime if $K$ has zero characteristic (see, e.g., \cite[Theorem~3.2.8]{JO}).
    \end{proof}
\end{thm}

\begin{rem}\label{remarkArch}
    Since $A$ consists of normal elements one obtains the following congruence relation $\cong$ on~$A$:
    \[a\cong b\iff\text{$a\mid nb$ and $b\mid ma$ for some positive integers $n$ and $m$.}\] The equivalence classes are
    subsemigroups of $A$; as in the commutative case we call these the Archimedean components of $A$. It follows
    that $A$ is a semilattice of its Archimedean components (the semilattice being $\Gamma=A/{\cong}$). From the
    proof of Theorem~\ref{semiprimeA} it follows that if $a\cong b$ then  $a,b\in A_e$ for some $e\in\Gamma$.
    Hence we obtain that each $A_e$ is a semilattice of the Archimedean components it contains. The theorem can
    thus be restated as: $K[A]$ is semiprime if and only if all Archimedean components of $A$ are cancellative.
    In the commutative case (i.e., if $(X,s)$ is the trivial solution) this is a well-known result
    (see, e.g., \cite[Corollary~9.6, Theorem~9.11 and Theorem~17.10]{Gilmer}).
\end{rem}

We link the decomposition as a semilattice of cancellative semigroups obtained in Theorem~\ref{semiprimeA}
with the chain obtained in Section~\ref{sec:desc} from divisibility. Note that $A=A(X,r)=M(X,s)$, where $s$
is the left  derived solution of $(X,r)$ and $a+b=b+\sigma_b (a)$, for some antihomomorphism $\sigma\colon(A,+)\to\Aut(A,+)$
(recall that indeed each $\sigma_a$ is bijective because $(X,r)$ is bijective). So the $\lambda$-map for the solution
$(X,s)$ is the constant map, i.e., $\lambda_x=\id$ for all $x\in X$. Again, as in the proof of Theorem~\ref{semiprimeA},
not to confuse with the multiplication in the semisimple algebra $Q=\Q_{cl}(K[A])$, we write the operation in $A$ also multiplicatively.
Hence with the notations as in Section~\ref{sec:desc}, if $\vn\ne Y\s X$ with $|Y|=i$ then either $A_{YY}=\vn$ or $\mc{L}_u=\vn$, i.e.,
$A_i/A_{i+1}$ is nil, or the equivalence classes of $\sim$ on $\mc{L}_u(i)$ are singletons, i.e., $A_i/A_{i+1}$
is an orthogonal $0$-disjoint union of semigroups $A_{YY}$ and its nil radical. Because $(X,s)$ is bijective, we know
(see, e.g., \cite{JKVA2018}) that for some $d\ge 1$ all elements $a^d$ for $a\in A$ are central. Let $c_Y=\prod_{y\in Y} y^d$,
a central element of $A$ that is divisible by all elements of $Y$. Consider the set $A_{YY}$ consisting of all elements of $A$
divisible by precisely the elements of $Y$ and assume $A_{YY}$ is a (non-empty) semigroup. Recall that for every $a\in A_{YY}$
there exists $b'\in A$ so that $ab'=c_Y^u$ for some $u\ge 1$. Also $a^u=c_Y b$ for some $b\in A$. Assume $c_Y\in A_{YY}$.
Hence also $a^u\in A_{YY}$. As $ a^u\in A_{YY}$ there exists $c\in A$ so that $a^uc=c_Y bc=c_Y^m$ for some $m\ge 1$.
Let $e$ be the central idempotent of $Q$ so that $Qc_Y=Qe$. Then $Qe=Qc_Y=Qc_Y^m=Qa^uQc=QaQc$.
So $Qe\s Qa$. Also $Qab'=Qc_Y^u=Qw_Y=Qe$. So $Qe\s Qa$. Hence $Qa=Qe$. Consequently, $A_{YY}\s A_e$.

On the other hand, assume $A_{YY}$ is not a semigroup and thus $A_{YY}$ is nil modulo $A_{i+1}$. Then, for $a\in A_{YY}$
there exists a positive integer, say $k$, so that $a^k\in A_{ZZ}$, for some non-empty subset $Z$ of $X$ so that $A_{ZZ}$
is cancellative semigroup. Let $f\in Q$ be the central idempotent so that $Qc_Z=Qf$. By the previous $a^k\in A_f$, i.e.,
$Qa^k=Qf$. Hence also $Qa=Qf$ and thus $a\in A_f$.

The above shows that, for $e\in\Gamma$, 
\[A_e=\bigcup_{\substack{\vn\ne Y\s X\text{ such that }Qc_Y=Qe\\ \text{and }A_{YY}\text{ is cancellative}}}A_{YY}\cup B_e\]
with $B_e=\{a\in A:a^u\in A_{YY},\,A_{YY}\text{ is cancellative and }A_{YY}\s A_e\}$.

Fix $e\in\Gamma$ and take $a\in A_f$ and $b\in A_g$ for some $f,g\in\Gamma$ satisfying $e\le f,g$. Choose $h\in\Gamma$
such that $\sigma_b(a)\in A_h$. Then $b\sigma_b(a)\in A_gA_h\s A_{gh}$, but also $b\sigma_b(a)=ab\in A_fA_g\s A_{fg}$.
Hence $e\le fg=gh\le h$. In particular, if \[A_e'=\bigcup_{f\in\Gamma\text{ such that }e\le f}A_f\] then the natural
extension $(A,s_A)$ of the left derived solution $(X,s)$ of $(X,r)$ to $A=A(X,r)=M(X,s)$ restricts to a solution on $A_e'$
and $(X,s)$ restricts to a solution on $X_e=A_e'\cap X$. Finally, the cancellative semigroup $A_e$ is an ideal in
the semigroup $A_e'=\free{B_e\cap X}$ and the structure group $G(X_e,s|_{X_e^2})\cong A_eA_e^{-1}$.

\section{The cancellative congruences and prime ideals}\label{sec:cong}

In the first part of this  section we describe the cancellative congruences of $A=A(X,r)$ and $M=M(X,r)$, for a finite left non-degenerate 
solution $(X,r)$ of the Yang--Baxter equation. Recall that the cancellative congruence of a semigroup $S$ is the smallest (with respect
to inclusion) congruence $\eta$ on $S$ such that $S/\eta$ is a cancellative semigroup. It turns out that the ideal of $K[M]$ determined
by this congruence is contained in every prime ideal $P$ of $K[M]$ with $P\cap M=\vn$. In the second part of this section we obtain crucial
information on all prime ideals of $K[A]$ which leads to an exact formula for the Gelfand--Kirillov dimension and classical Krull dimension
of $K[M]$ in a purely combinatorial way (in terms of  actions of certain finite monoids naturally related to the solution $(X,r)$).
Our results again extend the results obtained in \cite{JKVA2018} where $(X,r)$ also is assumed bijective.

In \cite{CJV21} the least left cancellative congruence on $A(X,r)$ and $M(X,r)$ has been recursively described
for left non-degenerate solutions $(X,r)$ and as an application one obtains the description of \cite{JKVA2018}
for finite bijective left non-degenerate solutions. The description obtained in this section is more explicit
(no iterative process is needed) and is based on a reduction to the cancellative congruence on $A(X,r)$. 

We first deal with $A=A(X,r)$. Let $v$ be the positive integer $v$ defined in \eqref{defv} (actually for our purposes
in this section we could work with the integer $e_{\sigma}$ defined in \eqref{idempotentexponent}). In particular, 
$\sigma_x^{v}=\sigma_{vx}$ is an idempotent mapping. Notice that if  $\sigma_a^2=\sigma_a$ and $\sigma_b^2=\sigma_b$,
for $a,b\in A$, then \[(\sigma_a\sigma_b)^2=\sigma_a\sigma_b\sigma_a\sigma_b=\sigma_{\sigma_a(b)}\sigma_a^2\sigma_b
=\sigma_{\sigma_a(b)}\sigma_a\sigma_b=\sigma_a\sigma_b^2=\sigma_a\sigma_b.\]
Hence it follows that $\sigma_a$ is idempotent for each $a\in\free{vx:x\in X}\s A$.

Recall from \eqref{defz} that \[z_{\id}=vx_1+\dotsb+vx_n\in A\]
and from \eqref{zcomm} that \[a+b+z_{\id}=b+a+z_{\id},\] for all $a,b\in\free{vx:x\in X}\s A$.

\begin{prop}\label{prop:7}
    Assume $(X,r)$ is a finite left non-degenerate solution of the Yang--Baxter equation. Let $A=A(X,r)$ and
    \[\eta_A=\{(a,b)\in A\times A:a+c=b+c \text{ for some }c\in A\}.\] Then $\eta_A$ is the cancellative congruence of $A$
    and the induced solution on $\ov{X}$, the image of $X$ in $\ov{A}=A/\eta_A$, is non-degenerate and bijective. In particular,
    $\ov{A}$ is a homomorphic image of the structure monoid of a non-degenerate bijective solution of the Yang--Baxter equation.
    Moreover, if $\kappa\in\Sym(n)$, where $n=|X|$ and \[\eta_i=\{(a,b)\in A\times A:a+iz_\kappa=b+iz_\kappa\}\] for $i\ge 1$,
    then $\eta_A=\bigcup_{i\ge 1}\eta_i$ and, as the monoid $A$ is left Noetherian, we have $\eta_A=\eta_t$ for some $t\ge 1$.
    \begin{proof}
        We first show that $\eta_A$ is a congruence. Let $(a,b)\in \eta_A$ and $w\in A$. Then, there exists $c\in A$
        such that $a+c=b+c$. Clearly, $w+a+c=w+b+c$, so $(w+a,w+b)\in \eta_A$. Moreover,
        \[a+w+c=a+c+\sigma_c(w)=b+c+\sigma_c(w)=b+w+c,\] which shows that $(a+w,b+w)\in \eta_A$ as well. 
        Next, we show that $\eta_A=\bigcup_{i\ge 1}\eta_i$. To do so, let $(a,b)\in \eta_A$.
        Then, there exists $c\in A$ such that $a+c=b+c$. Writing $X=\{x_1,\dotsc,x_n\}$, and ordering this set in an appropriate way,
        we may assume that $\kappa=\id$. Hence, let $z=z_\kappa=z_{\id}$. Then there exist non-negative integers $k_1,\dotsc,k_n$
        such that $c=k_1x_1+\dotsb+k_nx_n$. Let $tv$ be the smallest positive multiple of $v$ that is greater than
        $\max\{k_1,\dotsc,k_n\}$. Define
        \[c'=(tv-k_n)x_n+(tv-k_{n-1})\sigma_{tvx_n}(x_{n-1})+\dotsb+(tv-k_1)\sigma_{tvx_2+\dotsb+tvx_n}(x_1).\]
        Note that
        \begin{align*}
            c+c' & =w_1+k_n x_n+(tv-k_n)x_n+(tv -k_{n-1})\sigma_{tvx_n}(x_{n-1})+w_2\\
            & =w_1+tvx_n+(tv-k_{n-1})\sigma_{tvx_n}(x_{n-1})+w_2\\
            & =w_1+(vd-k_{n-1})x_{s-1}+tvx_n+w_2\\
            & =k_1x_1+\dotsb+k_{n-2}x_{n-2}+tdx_{n-1}+tdx_{n}+w_2
        \end{align*}
        with  $w_1=k_1x_1+\dotsb+k_{n-1}x_{n-1}$ and 
        \[w_2=(tv-k_{n-1})\sigma_{tvx_{n-1}+tvx_n}(x_{n-2})+\dotsb+(tv-k_1)\sigma_{tvx_2+\dotsb+tvx_n}(x_1).\]
        An induction argument yields that $c+c'=tvx_1+\dotsb+tvx_n$.  In particular, as $f+g+z=g+f+z$ for
        $f,g\in\free{vx:x\in X}$ by \eqref{zcomm}, one can reorder the first terms in $c+c'+z$. Hence we obtain
        \begin{align*}
            c+c'+z & =tvx_1+\dotsb+tvx_n+z\\
            & =vx_1+vx_2+\dotsb+vx_n+(t-1)vx_1+\dotsb+(t-1)vx_n+z=\dotsb=(t+1)z.
        \end{align*}
        Hence \[a+(t+1)z=a+c+c'+z=b+c+c'+z=b+(t+1)z.\] So, indeed, $\eta_A=\bigcup_{i\ge 1}\eta_i$.
        Clearly, $\eta_A$ is the least right cancellative congruence. We now show that it also is a left cancellative congruence.
        For this, suppose that $\ov{y}+\ov{a}=\ov{y}+\ov{b}$ for some $a,b,y\in A$ (here and later in the proof $\ov{c}$ denotes the image
        of $c\in A$ in $\ov{A}$). Then, in $A$ we have
        that there exists a positive integer $k$ such that $y+a+kz=y+b+kz$. We claim that we may replace $y$ by
        an element of $V=\free{vx:x\in X}\s A$. Indeed, let $y=t_1y_1+\dotsb+t_ny_n$ for some distinct $y_1,\dotsc,y_n\in X$ and some
        non-negative integers $t_1,\dotsc,t_n$. Then one can see inductively that $(v-t_n)y_n+\dotsb+(v-t_1)y_1+y\in V$.
        Hence we may assume that $y+a+kz=y+b+kz$  with $y\in V$. As, for any $w_1,w_2\in V$,
        we have $w_1+w_2+z=w_2+w_1+z$, it follows that there exists a positive integer $j$ and $u\in V$ 
        such that $jz=u+y+z$. Then
        \begin{align*}
            (j-1) z+a+kz &=j z+\sigma_z(a)+(k-1)z\\
            & =u+y+z+\sigma_z(a)+(k-1)z\\
            & =u+y+a+kz=u+y+b+kz\\
            & =u+y+b+z+(k-1)z\\
            & =u+y+z+\sigma_z(b)+(k-1)z\\
            & =j z+\sigma_z(b)+(k-1)z\\
            & =(j-1)z+z+\sigma_z(b)+(k-1)z\\
            & =(j-1)z+b+kz.
        \end{align*}
        Moreover, because $\sigma_z$ is idempotent and thus 
        \[a+kz=z+\sigma_z(a)+(k-1)z=z+\sigma_z^2(a)+(k-1)z=\sigma_z(a)+kz,\] 
        we obtain \[(j-1)z+a+kz=(j-1)z+\sigma_z(a)+kz=a+(j+k-1)z.\]
        And thus also \[(j-1)z+b+kz=(j-1)z+\sigma_z(b)+kz=b+(j+k-1)z.\]
        This shows that $(a,b)\in \eta_A$. Hence in $\ov{A}$, it holds that $\ov{a}=\ov{b}$.
        In particular, $\ov{A}$ is a cancellative monoid. 
        
        On $\ov{X}$, the image of $X$ in $\ov{A}$, one can define the solution $\ov{r}$ by $\ov{r}(\ov{x},\ov{y})=(\ov{y},\ov{\sigma_y(x)})$.
        Indeed, let $(x,x'),(y,y')\in \eta_A$ with $x,x',y,y'\in X$. Then there exist $i,j\ge 1$ such that $x+iz=x'+iz$
        and $y+jz=y'+jz$. Moreover, as $\sigma_{x'}(z)\in V$, there exists $w\in A$ and a positive
        integer $l$ such that $\sigma_{x'}(jz)+w+z=lz$. As $\sigma_{iz}$ is idempotent, it follows that
        \begin{align*}
            \sigma_x(y)+iz+lz &=iz+\sigma_{x+iz}(y)+lz\\
            & =iz+\sigma_{x'+iz}(y)+lz\\
            & =\sigma_{x'}(y)+iz+lz\\
            & =\sigma_{x'}(y)+lz+iz\\
            & =\sigma_{x'}(y)+\sigma_{x'}(jz)+w+z+iz\\
            & =\sigma_{x'}(y+jz)+w+z+iz\\
            & =\sigma_{x'}(y'+jz)+w+z+iz\\
            & =\sigma_{x'}(y')+iz+lz.
        \end{align*}
        Thus $(\sigma_x(y),\sigma_{x'}(y'))\in \eta_A$, as desired.
        
        Note that $\ov{A}$ is a homomorphic image of the structure monoid of this solution on $\ov{X}$.
        Furthermore, as $a+vx=vx+\sigma_{vx}(a)$ in $A$, it follows that $\sigma_{vx}$ induces a bijective
        map on the cancellative monoid $\ov{A}$ and hence on $\ov{X}$, showing that the solution $(\ov{X},\ov{r})$
        is non-degenerate and bijective.
    \end{proof}
\end{prop}

\begin{prop}\label{cancellativelnd}
    Assume $(X,r)$ is a left non-degenerate solution of the Yang--Baxter equation. If $\eta_A$ denotes the
    cancellative congruence of $A=A(X,r)$ and $M=M(X,r)=\{(a,\lambda_a):a\in A\}$ then
    \begin{align*}
        \eta_M &=\{(x,y)\in M\times M:x\circ z=y\circ z\text{ for some }z\in M\}\\
        & =\{((a,\lambda_a),(b,\lambda_b))\in M\times M:(a,b)\in\eta_A\text{ and }\lambda_a=\lambda_b\}
    \end{align*}
    is the cancellative congruence of $M$.
    \begin{proof}
        We consider $M$ as a submonoid of the semidirect product $A\rtimes\free{\lambda_x:x\in X}$.
        First, we shall show that the two sets of the statement are equal. Suppose
        $(a,\lambda_a)\circ(c,\lambda_c)=(b,\lambda_b)\circ(c,\lambda_c)$ for some $a,b,c\in A$. Then
        $\lambda_a\lambda_c=\lambda_b\lambda_c$ and thus $\lambda_a=\lambda_b$ because $\lambda_c$ is bijective.
        Furthermore, it follows that $a+\lambda_a(c)=b+\lambda_b(c)$ in $A$. Since $\lambda_a=\lambda_b$, one obtains that 
        $a+\lambda_a(c)=b+\lambda_a(c)$, which shows that $(a,b)\in\eta_A$. Conversely, suppose that $a,b\in A$ satisfy
        $(a,b)\in\eta_A$ and $\lambda_a=\lambda_b$. Then there exists $c\in A$ such that $a+c=b+c$. In particular, using
        $\lambda_a=\lambda_b$, we get
        \[(b,\lambda_b)\circ (\lambda_a^{-1}(c),\lambda_{\lambda_a^{-1}(c)})=(b+c,\lambda_{b+c})
        =(a+c,\lambda_{a+c})=(a,\lambda_a)\circ (\lambda_a^{-1}(c),\lambda_{\lambda_a^{-1}(c)}),\]
        which shows that the reverse inclusion holds as well.
        
        Clearly the sets yield a left congruence. We now show they also describe a right congruence, and thus
        describe the least cancellative congruence $\eta_M$. Abusing notation, we denote the mentioned sets
        already as $\eta_M$. To do so, let $((a,\lambda_a),(b,\lambda_b))\in \eta_M$ and $(c,\lambda_c)\in M$.
        Then, for some positive integer $i$, it holds that $a+iz=b+iz$, where $z=z_{\id}$ and $\lambda_a=\lambda_b$. We need to show that
        \[((a+\lambda_a(c),\lambda_a\lambda_c),(b+\lambda_b(c),\lambda_b\lambda_c))=((a,\lambda_a)\circ(c,\lambda_c),
        (b,\lambda_b)\circ(c,\lambda_c))\in\eta_M.\] As $\lambda_a=\lambda_b$, it follows that
        $\lambda_a\lambda_c=\lambda_b\lambda_c$. Moreover, $a+\lambda_a(c)=a+\lambda_b(c)$. This entails that
        \[a+\lambda_a(c)+iz=a+\lambda_b(c)+iz=a+iz+\sigma_{iz}(\lambda_b(c))=b+iz+\sigma_{iz}(\lambda_b(c))=b+\lambda_b(c)+iz.\]
        So, indeed, $(a+\lambda_a(c),b+\lambda_b(c))\in\eta_A$, and thus $\eta_M$ is the right cancellative congruence.
        
        Finally, we show that the quotient $\ov{M}=M/\eta_M$ is also left cancellative. Let
        $(a,\lambda_a),(b,\lambda_b),(c,\lambda_c)\in M$ be such that
        $\ov{(c,\lambda_c)}\circ\ov{(a,\lambda_a)}=\ov{(c,\lambda_c)}\circ\ov{(b,\lambda_b)}$ in $\ov{M}$
        (here and later in the proof $\ov{m}$ denotes the image of $m\in M$ in $\ov{M}$).
        Replacing $\ov{(c,\lambda_c)}$ by its certain power we may assume that $\lambda_c=\id$. Then, there exists
        $(w,\lambda_w)\in M$ such that \[(c,\id)\circ(a,\lambda_a)\circ(w,\lambda_w)=(c,\id)\circ(b,\lambda_b)\circ(w,\lambda_w).\]
        This is equivalent to $c+a+\lambda_a(w)=c+b+\lambda_b(w)$ and $\lambda_a=\lambda_b$. Due to the latter
        $\lambda_a(w)=\lambda_b(w)$, which entails that there exists $w'\in A$ such that $c+a+w'=c+b+w'$.
        As $\eta_A$ is the cancellative congruence of $A$, it follows that $(a,b)\in\eta_A$.
        In turn, this shows that $((a,\lambda_a),(b,\lambda_b))\in\eta_M$.
    \end{proof}
\end{prop}

In particular, the cancellative congruence of $M(X,r)$ can be made further explicit for non-degenerate bijective solutions $(X,r)$.
This sharpens \cite[Proposition~4.2]{JKVA2018}.

\begin{cor}\label{cor:cancon}
    Let $(X,r)$ be a non-degenerate bijective solution of the Yang--Baxter equation.
    If $\eta_A$ denotes the cancellative congruence of $A=A(X,r)$ then the cancellative congruence of $M=M(X,r)$ is given by
    \[\eta_M=\{((a,\lambda_a),(b,\lambda_b))\in M\times M:(a,b)\in \eta_A\}.\]
    \begin{proof}
        Let $(a,b)\in\eta_A$. By Proposition~\ref{prop:7} there exists a positive integer $i$ such that $a+iz=b+iz$, where $z=z_{\id}$.
        Moreover, since the element $z$ is contained in the set $\{m\in M:\lambda_m=\sigma_m=\id\}$, which is shown in \cite[Lemma~2]{MR4105532}
        to be $\lambda$-invariant because the solution is bijective by assumption, we get $\lambda_{\lambda_a^{-1}(z)}=\lambda_{\lambda_b^{-1}(z)}=\id$.
        Hence, \[\lambda_a=\lambda_a\lambda_{\lambda_a^{-1}(z)}=\lambda_{a+iz}=\lambda_{b+iz}=\lambda_b\lambda_{\lambda_b^{-1}(z)}=\lambda_b.\]
        Therefore, $((a,\lambda_a),(b,\lambda_b))\in\eta_M$ by Proposition~\ref{cancellativelnd}, which shows the result.
    \end{proof}
\end{cor}

The following corollary shows that the cancellative congruence $\eta_M$ of $M=M(X,r)$ is important to understand prime ideals
of $K[M]$ that do not intersect the monoid $M$ for  finite  left non-degenerate solutions $(X,r)$ of the Yang--Baxter equation.
Also recall that since $K[M]$ is a PI-algebra, the Jacobson radical $J(K[M])$ equals its prime radical $B(K[M])$, see for example
\cite[Theorem~21.12]{Okn1991}. Furthermore this is a nilpotent ideal as $K[M]$ is an affine algebra. The left annihilator of
a subset $S$ of a ring $R$ is denoted $\lann_R(S)$.

\begin{cor}
    Assume $(X,r)$ is a left non-degenerate solution of the Yang--Baxter. Let $M=M(X,r)$. If $K$ is a field and $I(\eta_M)$
    denotes the ideal of $K[M]$ determined by the cancellative congruence $\eta_M$ of $M$ (i.e., the ideal generated by all
    elements $s-t$ with $(s,t)\in\eta_M$, or equivalently the $K$-vector space spanned by such elements) then:
    \begin{enumerate}
        \item $I(\eta_M)\s P$ for each prime ideal $P$ of $K[M]$ satisfying $P\cap M=\vn$,
        \item $I(\eta_M)=\lann_{K[M]}(M_{XX}^t)$, for some positive integer $t$.
    \end{enumerate}
    If $\ch K=0$ then:
    \begin{enumerate}
       \setcounter{enumi}{2}
        \item $I(\eta_M)=\bigcap P$, where the intersection runs over all prime ideals $P$ of $K[M]$ that intersect $M$ trivially, and thus
        \item $J(K[M])=B(K[M])=I(\eta_M)\cap\bigcap_{P\in\mc{P}}P$, where \[\mc{P}=\{P\in\Spec(K[M]):P\cap M\ne\vn\}=\{P\in\Spec(K[M]):M_{XX}\s P\}.\]
    \end{enumerate}
    \begin{proof}
        (1) Let $((a,\lambda_a),(b,\lambda_b))\in\eta_M$. Then, by Proposition~\ref{prop:7}, there exists a positive integer $i$ such that 
        $a+iz_\kappa=b+iz_\kappa$ for each $\kappa\in\Sym(n)$, where $n=|X|$. Let $(c,\lambda_c)\in M$. Write $X=\{x_1,\dotsc,x_n\}$ and
        put $z=z_{\id}$. Because $\lambda_a=\lambda_b$, by Proposition~\ref{cancellativelnd}, there exists a permutation $\kappa\in\Sym(n)$
        such that $\lambda_a(\lambda_c(x_i))=x_{\kappa(i)}=\lambda_b(\lambda_c(x_i))$ for each $1\le i\le n$. Then
        \begin{align*}
            (a,\lambda_a)\circ(c,\lambda_c)\circ(iz,\lambda_{iz}) &
            =(a+\lambda_a(c)+\lambda_a(\lambda_c(iz)),\lambda_a\lambda_c\lambda_{iz})\\
            & =(a+iz_\kappa+\sigma_{iz_\kappa}(\lambda_a(c)),\lambda_a\lambda_c\lambda_{iz})\\
            & =(b+iz_\kappa+\sigma_{iz_\kappa}(\lambda_b(c)),\lambda_b\lambda_c\lambda_{iz})\\
            & =(b+\lambda_b(c)+\lambda_b(\lambda_c(iz)),\lambda_b\lambda_c\lambda_{iz})\\
            & =(b,\lambda_b)\circ(c,\lambda_c)\circ(iz,\lambda_{iz}).
        \end{align*}
        Therefore, we conclude that $((a,\lambda_a)-(b,\lambda_b))\circ K[M]\circ(iz,\lambda_{iz})=0$. Since, by assumption,
        $(iz,\lambda_{iz})\notin P$, we get $(a,\lambda_a)-(b,\lambda_b)\in P$. Hence, $I(\eta_M)\s P$, which shows the first part of the result.

        (2) From the proof of (1) it follows that there exists a positive integer $t$ such that  $((a,\lambda_a),(b,\lambda_b))\in\eta_M$ if and only if 
        $((a,\lambda_a)-(b,\lambda_b))\circ(tz,\lambda_{tz})=0$ (and this for any $z=z_{\kappa}$ with $\kappa \in \Sym (n)$), or equivalently
        $((a,\lambda_a)-(b,\lambda_b))\circ(tz,\lambda_{tz})\circ M=0$. Clearly $(t^2z,\lambda_{t^2z})\circ M\s M_{XX}^{t^2}\s(tz,\lambda_{tz})\circ M$.
        It follows that $I(\eta_M)=\lann_{K[M]}(M_{XX}^{t})$.

        Clearly $K[M]/I(\eta_M) \cong K[T]$, where $T=M/\eta_M$ is a cancellative monoid. As an epimorphic image of a PI-algebra, the monoid algebra
        $K[T]$ also is PI. It is then  well-known (see, e.g., \cite[Theorem~3.1.9]{JO}) that $T$ has a group of fractions, say $H$, and $K[H]$ is
        a PI-algebra. As $M$ is finitely generated as a monoid, $H$ is finitely generated as a group and thus $H$ is abelian-by-finite group
        (see, e.g., \cite[Corollary 13]{Okn1991} or \cite[Theorem~5.3.15]{MR798076}). In particular, $K[H]$ is a Noetherian algebra and thus
        $B(K[T])=B(K[H])\cap K[T]$ (see for example \cite[Corollary 11.5]{Okn1991}).

        To prove (3) and (4), assume $\ch K=0$. Then $B(K[H])=0$ (see for example \cite[Theorem~3.2.8]{JO}), and thus $\{0\}$ is the 
        intersection of all prime ideals $P$ of $K[H]$. Furthermore, since $K[H]$ is Noetherian, it follows from \cite[Lemma~7.15]{Okn1991}
        that $P\cap K[T]$ is a prime ideal of $K[T]$ and clearly it does not intersect $T$. Consequently,  $\{0\}$ is the intersection of
        all prime ideals $Q$ of $K[T]$ with $T\cap Q= \vn$. Clearly such prime ideals $Q$ correspond to prime ideals of $K[M]$ that do not
        intersect $M$ and that contain $I(\eta_M )$. This proves (3). Part (4) is now clear as well.
    \end{proof}
\end{cor}

We now give a description of prime ideals of the monoid $A=A(X,r)$ and its algebra $K[A]$ over a field $K$. As the proofs
of analogous results for bijective solutions presented in \cite{JKVA2018} remain valid, we shall only formulate these
results omitting their proofs.

The prime spectra of $A$ and $K[A]$ are denoted as $\Spec(A)$ and $\Spec(K[A])$, respectively.

The authors were informed by Be'eri Greenfeld \cite{Greenfeld} that there is a surprising and interesting connection
between structure algebras of left derived non-degenerate (and thus bijective) solutions of the Yang--Baxter equation,
an extensively studied arithmetic-geometric conjecture (see \cite{EllWest} and \cite{MR3488737}) and a conjecture of
Malle \cite{MR2068887} that predicts the asymptotic behaviour of the number $N_G(K,D)$ of extensions of a number field
$K$  whose Galois closure has  Galois group  $G$ with the absolute value of the discriminant less than  $D$. This
conjecture has been extended to arbitrary global fields since then, and Ellenberg, Tran and Westerland in \cite{EllWest}
tackled the ``geometric Malle conjecture'', namely for $K=k(x)$ a function field over a finite field. They provide an upper
bound to $N_{G}(K,X)$ and   a suitable poly-logarithmic factor in their bound is obtained  via the Gelfand--Kirillov dimension
of a certain non-commutative graded algebra. Greenfeld \cite{Greenfeld} pointed out that these algebras can be realized as
structure algebras of non-degenerate (and thus bijective) left derived solutions of the Yang--Baxter equation. The
Gelfand--Kirillov dimension of these algebras has been calculated in \cite[Theorem~3.8]{JKVA2018}.
We now extend this result to left non-degenerate solutions that are not necessarily bijective. 

\begin{prop}\label{prop:6'}
	Assume $(X,r)$ is a left non-degenerate solution of the Yang--Baxter equation. Let $A=A(X,r)$ and
	\[\mc{Z}=\mc{Z}(X,r)=\{Z\s X:\vn\ne Z\ne X,\,\sigma_x(Z)\s Z\text{ and }
	\sigma_x(X\setminus Z)\s X\setminus Z\text{ for all }x\in X\setminus Z\}.\]
	For $Z\in\mc{Z}$ put $P(Z)=\bigcup_{z\in Z}(z+A)$. The maps
	\[\mc{Z}\to\Spec(A)\colon Z\mapsto P(Z)\quad\text{and}\quad\Spec(A)\to\mc{Z}\colon P\mapsto X\cap P\]
	are mutually inverse bijections.
\end{prop}

\begin{prop}\label{prop:8.3}
	Assume $(X,r)$ is a finite left non-degenerate solution of the Yang--Baxter equation. Let $A=A(X,r)$ and
	$\mc{Z}=\mc{Z}(X,r)$. If $K$ is a field and $P$ is a prime ideal of the algebra $K[A]$ then
	$X\cap P\in\mc{Z}\cup\{\vn,X\}$. Moreover, for such a prime ideal $P$ the following properties hold:
	\begin{enumerate}
		\item there exists an inclusion preserving bijection between the set of prime ideals $Q$ of $K[A]$ with
		the property $X\cap Q=X\cap P$ and the set of all prime ideals of the algebra $K[A\setminus P]$.
		Moreover, the monoid $A\setminus P$ has the following presentation
		\[A\setminus P\cong\free{X\setminus P\mid x+y=y+\sigma_y(x)\text{ for all }x,y\in X\setminus P},\]
		and thus it is the structure monoid of the subsolution of $(X,r)$ on $X\setminus P$.
		\item there exists an inclusion preserving bijection between the set of prime ideals $Q$ of $K[A]$ satisfying
		$Q\cap A=\vn$ and the set of all prime ideals of the group algebra $K[G]$, where
		\[G=\gr(X\mid x+y=y+\sigma_y(x)\text{ for all }x,y\in X).\]
		Furthermore, the cancellative monoid $\ov{A}=A/\eta_A$ has a group of quotients, isomorphic to $G$,
        which is a central localization of $\ov{A}$ with respect to a submonoid generated by a single element
        that is divisible by all generators. Clearly, $G$ is a finite conjugacy group.
	\end{enumerate}
\end{prop}

\begin{rem}
    Let $(X,r)$ be a finite left non-degenerate solution of the Yang--Baxter equation. Then $R=K[M(X,r)]$
    is an affine PI-algebra by \eqref{S1}. Hence, by a well-known result (see, e.g., \cite[Theorem~13.10.3]{MR}),
    $R$ is a Jacobson ring (i.e., every prime image of $R$ is semiprimitive), every simple (say) right $R$-module $V$
    is finite dimensional over $K$ and the algebra $R$ satisfies the Nullstellensatz (i.e., the endomorphism
    algebra $\End(V_R)$ is algebraic over $K$).
\end{rem}

In order to compute the Gelfand--Kirillov dimension of $K[M(X,r)]$ we need one more lemma that describes the
largest bijective non-degenerate homomorphic image of a left derived solution. 

\begin{lem}\label{lem:bijectivehom}
    Let $(X,r)$ be a left non-degenerate set-theoretic solution of the Yang--Baxter equation. Let $(X,s)$ denote the left derived solution of
    $(X,r)$. Then the relation \[x\sim y\iff\sigma_z(x)=\sigma_z(y),\] where $z=z_{\id}\in A=A(X,r)$, is an equivalence relation on $X$.
    Furthermore, the solution $(X,s)$ induces a bijective non-degenerate solution on the set $\ov{X}=X/{\sim}$. 
    \begin{proof}
        It is clear that $\sim$ is an equivalence relation on $X$. Let $x,x',y,y' \in  X$ with $x\sim x'$ and $y\sim y'$. Then
        \[\sigma_z(\sigma_x(y))=\sigma_{x+z}(y)=\sigma_{z+\sigma_z(x)}(y)=\sigma_{\sigma_z(x)}(\sigma_z(y))
        =\sigma_{\sigma_z(x')}(\sigma_z(y'))= \sigma_z(\sigma_{x'}(y')).\]
        Thus, $(X,s)$ induces a solution on the set $\ov{X}=X/{\sim}$. Moreover, for any $x\in X$ there exists $a\in A$ such that $z=x+a$.
        Hence, the map induced on $\ov{X}$ by $\sigma_x$ is invertible as $\sigma_z$ induces the identity on $\ov{X}$. Indeed, as
        $\sigma_z^2(x)=\sigma_z(x)$, it follows that $x\sim \sigma_z(x)$. This shows the result.
    \end{proof}
\end{lem}

The following result extends \cite[Theorem 3.5]{JKVA2018} and determines the classical Krull dimension of $K[A]/P$ for a minimal prime ideal $P$
of $K[A]$ not intersecting $A=A(X,r)$. To explicitize this number, we recall the notion of an orbit under the action of a monoid $S$ on a set $X$
(the image of $x\in X$ under the action of $s\in S$ is denoted by $s\cdot x$). We say $x,y\in X$ are in the same orbit if and only if there exist
$s,t\in S$ such that $s\cdot x=t\cdot y$. Moreover, if for any $s,s'\in S$ there exists $t\in S$ such that $ss'=ts$ then the relation of being
in the same orbit is an equivalence relation. The latter holds for the monoid $\Sigma=\free{\sigma_x:x\in X}$ because
$\sigma_x\sigma_y=\sigma_{\sigma_x(y)}\sigma_x$ for all $x,y\in X$.

\begin{prop}\label{prop:6.7}
    Assume $(X,r)$ is a finite left non-degenerate solution of the Yang--Baxter equation. Let $A=A(X,r)$. If $K$ is a field and $P$
    is a minimal prime ideal of $K[A]$ satisfying $P\cap A=\vn$ then \[\clK K[A]/P=\clK K[G]=t,\] where $G$ is as in Proposition~\ref{prop:8.3}
    and $t$ is the number of orbits of $X$ with respect to the action of the monoid $\Sigma=\free{\sigma_x:x\in X}$.
    \begin{proof}
        To avoid confusion of the binary operation of $A$ with multiplication in the algebra $K[A]$, the operation in $A$ will
        be denoted multiplicatively. As mentioned in previous sections, there exists a positive integer $d$ such that the map
        $\sigma_a^d$ is idempotent for each $a\in A$. First, we claim that for any $a,b\in A$ we have \[b-\sigma_a^d(b)\in P.\]
        Indeed, as $ba^d=a^d\sigma_a^d(b)=a^d\sigma_a^{2d}(b)=\sigma_a^d(b)a^d$, we get $(b-\sigma_a^d(b))K[A]a^d\s(b-\sigma_a^d(b))a^dK[A]=0$
        because $Aa^d\s a^dA$. Since, by assumption, $a^d\notin P$, the claim follows. In particular, using the relation $\sim$ defined in 
        Lemma~\ref{lem:bijectivehom}, one obtains that if $x\sim y$ for some $x,y\in X$ then $\sigma_z(x)=\sigma_z(y)$ and thus
        \[x-y=(x-\sigma_z^d(x))-(y-\sigma_z^d(y))\in P.\] Moreover, denoting the left derived solution of $(X,r)$ by $(X,s)$,
        the latter induces a finite bijective non-degenerate solution $(\ov{X},\ov{s})$ on the set $\ov{X}=X/{\sim}$. Let $\ov{P}$ denote
        the natural image of the ideal $P$ in the algebra $K[\ov{A}]$, where $\ov{A}=A(\ov{X},\ov{s})$. Then $K[A]/P\cong K[\ov{A}]/\ov{P}$.
        Hence, by \cite[Theorem 3.5]{JKVA2018}, one has that $$ \clK K[A]/P =\clK K[\ov{A}]/\ov{P}=t,$$ where $t$ is the number of orbits of
        $\ov{X}$ under the action of the group $\ov{\Sigma}$ generated by the permutations of $\ov{X}$ induced by $\sigma_x$ with $x\in X$.
        Clearly, as $\sim$ identifies elements inside the same orbit, the number of orbits of $X$ under the action of the monoid $\Sigma$
        is equal to $t$. Hence the result follows.
    \end{proof}
\end{prop}

Motivated by Proposition~\ref{prop:6.7} we define a monoid \[\Sigma_Z=\free{\sigma_x:x\in X\setminus Z}\]
and \[t(Z)=\text{the number of orbits of }X\setminus Z\text{ with respect to the action of }\Sigma_Z\]
for each $Z\in\mc{Z}_0=\mc{Z}\cup\{\vn\}$, where $\mc{Z}=\mc{Z}(X,r)$ (see Propositions~\ref{prop:6'} and \ref{prop:8.3}).

We conclude this section with a combinatorial description of the Gelfand--Kirillov dimension of $K[M(X,r)]$. This extends
\cite[Theorem~3.8]{JKVA2018}. Because of the results stated above the proof remains the same. However, since it is rather
short we include it for completeness' sake.

\begin{thm}\label{GKdim}
	Assume that $(X,r)$ is a finite left non-degenerate solution of the Yang--Baxter equation.
	Let $A=A(X,r)$, $M=M(X,r)$ and $\mc{Z}_0=\mc{Z}\cup\{\vn\}$, where $\mc{Z}=\mc{Z}(X,r)$. If $K$ is a field then
	\[\GK K[M]=\clK K[M]=\GK K[A]=\clK K[A]=\max\{t(Z):Z\in\mc{Z}_0\}\le|X|.\]
	\begin{proof}
        As before, we shall use the multiplicative notation for the binary operation in $A$. Because the finitely generated $K$-algebra
        $K[A]$ is left Noetherian by \eqref{N1} and PI by \eqref{S1}, \cite[Theorem 3.5.2]{JO} implies that $\GK K[A]=\clK K[A]$.
        So, it remains to show that $\clK K[A]=t$, where \[t=\max\{t(Z):Z\in\mc{Z}_0\}.\] If $P$ is a minimal prime ideal of $K[A]$
        then $X\cap P\in\mc{Z}_0$ and $P$ corresponds to a minimal prime ideal $Q$ of the algebra $K[A\setminus P]$ such that
		$Q\cap(A\setminus P)=\vn$ and $K[A]/P\cong K[A\setminus P]/Q$ (see Proposition~\ref{prop:8.3}).
		Therefore, Proposition~\ref{prop:6.7} yields \[\clK K[A]/P=\clK K[A\setminus P]/Q=t(X\cap P)\le t.\]
		Since we have $\clK K[A]=\clK K[A]/P$ for some minimal prime ideal $P$ of $K[A]$, the inequality $\clK K[A]\le t$ follows.
		To show that $\clK K[A]\ge t$ we have to check that $\clK K[A]\ge t(Z)$ for each $Z\in\mc{Z}_0$.
		So, let us fix $Z\in\mc{Z}_0$. If $Z=\vn$ then we are done by Proposition~\ref{prop:6.7}.
		Whereas, if $Z\in\mc{Z}$ then $A=P(Z)\cup A(Z)$, where $P(Z)=\bigcup_{z\in Z}zA$ and
		$A(Z)=A\setminus P(Z)=\free{X\setminus Z}$ is the submonoid of $A$ generated by $X\setminus Z$.
		Therefore, $K[A]/K[P(Z)]\cong K[A(Z)]$ and this leads to
		\[\clK K[A]\ge\clK K[A]/K[P(Z)]=\clK K[A(Z)]\ge t(Z).\]
		The last inequality is a consequence of the fact that the ideal $P_0$ of $K[A(Z)]$, generated
		by elements of the form $x-y$ for all $x,y\in X\setminus Z$ which are in the same orbit of $X\setminus Z$
		with respect to the action of the monoid $\Sigma_Z$, satisfies $K[A(Z)]/P_0\cong K[x_1,\dotsc,x_{t(Z)}]$,
		the polynomial algebra in $t(Z)$ commuting variables.
  
        Clearly $\GK K[M]=\GK K[A]$, as the bijective $1$-cocycle $\pi\colon M\to A$ is degree preserving. Further,
        $K[M]$ is left Noetherian by Theorem~\ref{theorem:A} and thus $\GK K[M]=\clK K[M]$ by \cite[Theorem~3.5.2]{JO}.
        Hence the result follows.
	\end{proof}
\end{thm}

As a consequence of Theorems~\ref{primefinitegroup} and \ref{GKdim} we obtain a characterization of involutive solutions
in the class of all finite left non-degenerate solutions with two-sided Noetherian structure algebras. As this class
contains all solutions which are additionally bijective, the following result is a substantial strengthening  
of results obtained in \cite[Theorem~4.6]{JKVA2019} (see also \cite[Theorem 4.5]{JKVA2018}).

Recall that the rank $\rk S$ of a monoid $S$ is the largest possible rank of a free abelian submonoid of $S$ (see \cite{CP}
or \cite{Okn98,Okn1991,JO} for details). For definitions of all homological notions (including the global dimension $\gld R$ and
the left and right injective dimension of a two-sided Noetherian ring $R$, in case both the latter dimensions coincide we denote
this number by $\id R$) used below we refer to \cite{BG,MR4278764,Lav,SZ}.

\begin{cor}\label{cor:primefinitegroup}
    Let $(X,r)$ be a finite left non-degenerate solution of the Yang--Baxter equation.
    If $K$ is a field and $M=M(X,r)$ then the following conditions are equivalent:
    \begin{enumerate}
        \item $(X,r)$ is an involutive solution.
        \item $M$ is a cancellative monoid and $\Omega_\lambda$ is a trivial group.
        \item $K[M]$ is a prime algebra and $\Omega_\lambda$ is a trivial group.
        \item $K[M]$ is a domain.
        \item $\GK K[M]=|X|$.
    \end{enumerate}
    Moreover, if the algebra $K[M]$ is right Noetherian, then the above conditions are equivalent to:
    \begin{enumerate}
        \setcounter{enumi}{5}
        \item $\rk M=|X|$.
        \item $\clK K[M]=|X|$.
        \item $\id K[M]=|X|$.
        \item $K[M]$ has finite global dimension.
        \item $K[M]$ is an Auslander--Gorenstein algebra with $\id K[M]=|X|$.
        \item $K[M]$ is an Auslander-regular algebra.
    \end{enumerate}
    \begin{proof}
        The equivalence of conditions (1)--(4) was already proved in Theorem~\ref{primefinitegroup}.
        Next, the implication $(1)\Longrightarrow(5)$ easily follows by \cite[Theorem 1.6]{GVdB}.
        If $\GK K[M]=|X|$ then Theorem~\ref{GKdim} assures that $X$ decomposes into $|X|$ orbits under the action
        of the monoid $\Sigma=\free{\sigma_x:x\in X}$. Therefore, the orbits are singletons and it follows that
        $\sigma_x=\id$ for each $x\in X$. Hence $r^2=\id$ by \eqref{R3}, and thus we get $(5)\Longrightarrow(1)$.
        Summarizing, we have already shown that all conditions (1)--(5) are equivalent.

        From now on, we assume that the algebra $K[M]$ also is right Noetherian. Because $K[M]$ is also representable
        by \eqref{S1}, it follows that $M$ is a linear monoid. Hence Theorem~\ref{GKdim} and
        \cite[Proposition~1, p.~221, Proposition~7, p.~280–281, and Theorem~14, p.~284]{Okn1991} yield
        \[\GK K[M]=\clK K[M]=\rk M.\] In particular, conditions (5), (6) and (7) are equivalent.
        So, in consequence, all statements (1)--(7) are equivalent. Next, it is well-known that
        $(1)\Longrightarrow(9)$ (see \cite[Theorem~1.4]{GVdB}). Further, if $\id K[M]<\infty$ then
        \cite[Theorem~1, p.~126]{BG} implies that the algebra $K[M]$ is Auslander--Gorenstein and
        \begin{equation}
            \GK K[M]=\clK K[M]=\id K[M].\label{homdim}
        \end{equation}
        Moreover, if $\gld K[M]<\infty$ then \cite[Theorem~1, p.~126]{BG} says also that $K[M]$ is an Auslander-regular domain.
        In particular, we get $(9)\Longrightarrow(11)$. Further, $(11)\Longrightarrow(10)$ as well, because in this case $K[M]$
        is a domain, so by the implication $(4)\Longrightarrow(5)$ and \eqref{homdim} we obtain $\id K[M]=\GK K[M]=|X|$.
        Clearly, $(10)\Longrightarrow(8)$. Finally, $(8)\Longrightarrow(5)$ again by \cite[Theorem~1, p.~126]{BG}.
        Hence the result is proved.
	\end{proof}
\end{cor}

\section{Structure algebras of some degenerate solutions}\label{sec:allrhoequal}

So far we focused on finite solutions that are left non-degenerate. In this section we make a first attempt to deal with
degenerate solutions. We consider arbitrary finite solutions $(X,r)$ of the Yang--Baxter equation with
\begin{equation}\label{fixedrho}
    r(x,y)=(\lambda_x(y),\rho(x)),
\end{equation}
i.e., all maps $\rho_x$ are equal to a fixed map $\rho$. Note that this means 
\[\lambda_x\lambda_y=\lambda_{\lambda_x(y)}\lambda_{\rho(x)}\quad\text{and}\quad\rho\lambda_x=\lambda_{\rho(x)}\rho\]
for all $x,y\in X$. Examples of such left non-degenerate solutions are solutions associated to metahomomorphism on groups.
The definition of metahomomorphism for a group has been introduced in \cite{MR1613779}. Here we recall the definition.
Let $X$ be a (finite) group and $\rho\colon X\to X$ an arbitrary function. Then the map $r\colon X\times X\to X\times X$,
defined as $r(x,y)=(xy\rho(x)^{-1},\rho(x))$, is a left non-degenerate solution of the Yang--Baxter equation if and only if
$\rho$ satisfies \[\rho(xy\rho(x)^{-1})=\rho(x)\rho(y)\rho^2(x)^{-1}\] for all $x,y\in X$. Such a function $\rho$ is called
a \emph{metahomomorphism}. Note that for a solution $(X,r)$ of the above form we have $\lambda_x(y)=xy\rho(x)^{-1}$ and thus
$\lambda_x^{-1}(y)=x^{-1}y\rho(x)$. In particular, $\lambda_x^{-1}(x)=\rho(x)$ and
$\Lambda=\{q(x)=\lambda_x^{-1}(x):x\in X\}=\Img\rho$. Note also that
\[r(\rho(x),\rho(y))=(\rho(x)\rho(y)\rho^2(x)^{-1},\rho^2(x))=(\rho(xy\rho(x)^{-1}),\rho^2(x)),\]
and thus $r$ restricts to $\Lambda$. Examples of metahomomorphims are constant maps, endomorphisms and the inversion map
$x\mapsto x^{-1}$. There are many more examples, even for abelian groups (see \cite{MR2325536}). Another class of examples 
of the type \eqref{fixedrho} are the Lyubashenko solutions (see \cite{Dri1992}). These are solutions $(X,r)$ with $X$ an
arbitrary finite set and $r(x,y)=(\lambda(y),\rho(x))$, where $\lambda\colon X\to X$ and $\rho\colon X\to X$ are commuting
functions. Such a solution is left non-degenerate when $\lambda$ is a permutation. 

\begin{prop}\label{7.11}
    Let $(X,r)$ be a finite solution of the Yang--Baxter equation of the form $r(x,y)=(\lambda_x(y),\rho(x))$.
    If $K$ is a field and $M=M(X,r)$ then $K[M]$ is a right Noetherian PI-algebra of finite Gelfand--Kirillov dimension.
    \begin{proof}
        Since $X$ is a finite set, there exists $n\ge 1$ such that $\Img\rho^n=\Img\rho^{n+1}$; let us denote this set by $Y$.
        Since $(X,r)$ is a solution of the Yang--Baxter equation and $\rho_x=\rho$ for all $x\in X$, we know that
        $\lambda_{\rho(x)}\rho=\rho\lambda_x$. Hence, by induction, it follows that 
        \begin{equation}
            \lambda_{\rho^i(x)}\rho^i=\rho^i\lambda_x\quad\text{for each }i\ge 0\text{ and }x\in X.\label{7.11-1}
        \end{equation}
        In particular, \eqref{7.11-1} assures that $\lambda_y(Y)\s Y$ for each $y\in Y$, and thus $(X,r)$ restricts
        to a solution $(Y,r_Y)$ on $Y$. Moreover, the restriction $\rho_Y\colon Y\to Y$ of $\rho$ is surjective.
        Hence $\rho_Y$ is a bijection and thus $(Y,r_Y)$ is a right non-degenerate solution. Moreover, the diagonal
        map $q\colon Y\to Y\colon y\mapsto\rho_Y^{-1}(y)$ of the solution $(Y,r_Y)$ is bijective. Hence by the
        left-right dual of \eqref{N2}, $R=K[M(Y,r_Y)]$ is a right Noetherian PI-algebra of finite Gelfand--Kirillov dimension. 
        
        We claim that $K[M]$ is a finite right $R$-module, which will imply that $K[M]$ is right Noetherian PI-algebra.
        To prove this, first, note that for each choice of $1\le j<k$ and $x_1,\dotsc,x_k\in X$  there
        exist $x_1',\dotsc,x_k'\in X$ such that
        \begin{equation}\label{7.11-2}
            x_1\circ\dotsb\circ x_k=x_1'\circ\rho(x_2')\circ\dotsb\circ\rho^{j-1}(x_j')\circ\rho^j(x_{j+1}')\circ\dotsb\circ\rho^j(x_k').
        \end{equation}
        (Counting from the left, on the right-hand side of \eqref{7.11-2}, we have $j$ elements that are images of consecutive
        increasing powers of $\rho$ up to $\rho^{j-1}$, whereas all other elements on the right are in $\Img\rho^j$.)
        To prove \eqref{7.11-2} define \[S=\{(j,k)\in\mb{N}\times\mb{N}:j<k\}\] and
        \[(j_1,k_1)<(j_2,k_2)\iff k_1<k_2\text{ or }k_1=k_2\text{ but }j_1<j_2\] for $(j_1,k_1),(j_2,k_2)\in S$
        (we write $s\le t$ for $s,t\in S$ in case $s<t$ or $s=t$). Clearly $(S,\le)$ is a well-ordered set.
        Note that \eqref{7.11-2} holds for the smallest element $(j,k)=(1,2)$ of $S$ because $x_1\circ x_2=\lambda_{x_1}(x_2)\circ\rho(x_1)$.
        Now, fix $(1,2)<(j,k)\in S$ and assume that \eqref{7.11-2} holds for each $(j',k')\in S$ satisfying $(j',k')<(j,k)$.
        There are two cases to consider.
        
        (1) If $j=1$ then using $(k-1)$-times the fact that for all $x,x'\in X$ there exists $y\in X$ such that
        $x\circ x'=y\circ\rho(x)$ (clearly it is enough to put $y=\lambda_x(x')$), we see that there exist
        $y_1,\dotsc,y_{k-1}\in X$ such that
        \begin{align*}
            x_1\circ\dotsb\circ x_{k-1}\circ x_k & =x_1\circ\dotsb\circ x_{k-2}\circ y_1\circ\rho(x_{k-1})\\
            & =x_1\circ\dotsb\circ x_{k-3}\circ y_2\circ\rho(x_{k-2})\circ\rho(x_{k-1})\\
            & \vdotswithin{=}\\
            & =x_1\circ y_{k-2}\circ\rho(x_2)\circ\dotsb\circ\rho(x_{k-2})\circ\rho(x_{k-1})\\
            & =y_{k-1}\circ\rho(x_1)\circ\rho(x_2)\circ\dotsb\circ\rho(x_{k-2})\circ\rho(x_{k-1}).
        \end{align*}
        Hence \eqref{7.11-2} holds in case (1) with $x_1'=y_{k-1}$ and $x_i'=x_{i-1}$ for $2\le i\le k$.
        
        (2) If $j>1$ then $(j-1,k)\in S$ and $(j-1,k)<(j,k)$. So, by our assumption, \eqref{7.11-2} holds for $(j-1,k)$.
        Hence
        \begin{equation}\label{7.11-3}
            x_1\circ\dotsb\circ x_k=y_1\circ\rho(y_2)\circ\dotsb\circ\rho^{j-2}(y_{j-1})\circ\rho^{j-1}(y_j)\circ\dotsb\circ\rho^{j-1}(y_k)
        \end{equation}
        for some $y_1,\dotsc,y_k\in X$. Now, observe that if $x,y\in X$ and $0\le p\le q$ then
        \begin{align*}
            \rho^p(x)\circ\rho^q(y) & =\lambda_{\rho^p(x)}(\rho^p(\rho^{q-p}(y)))\circ\rho^{p+1}(x)\\
            & =\rho^p(\lambda_x(\rho^{q-p}(y)))\circ\rho^{p+1}(x).
        \end{align*}
        Applying the above to the right-hand side of \eqref{7.11-3} we get
        \begin{align*}
            y_1\circ{\rho(y_2)}\circ\dotsb & \circ\rho^{j-2}(y_{j-1})\circ\rho^{j-1}(y_j)\circ\dotsb\circ\rho^{j-1}(y_{k-1})\circ\rho^{j-1}(y_k)\\
            & =y_1\circ{\rho(y_2)}\circ\dotsb\circ\rho^{j-2}(y_{j-1})\circ\rho^{j-1}(y_j)\circ\dotsb\circ\rho^{j-1}(z_1)\circ\rho^j(y_{k-1})\\
            & \vdotswithin{=}\\
            & =y_1\circ{\rho(y_2)}\circ\dotsb\circ\rho^{j-2}(y_{j-1})\circ\rho^{j-1}(y_j)
            \circ\rho^{j-1}(z_{k-j-1})\circ\rho^j(y_{j+1})\circ\dotsb\circ\rho^j(y_{k-1})\\
            & =y_1\circ{\rho(y_2)}\circ\dotsb\circ\rho^{j-2}(y_{j-1})\circ\rho^{j-1}(z_{k-j})
            \circ\rho^j(y_j)\circ\rho^j(y_{j+1})\circ\dotsb\circ\rho^j(y_{k-1})
        \end{align*}
        for some $z_1,\dotsc,z_{k-j}\in X$. Therefore, \eqref{7.11-2} also holds in case (2) with $x_i'=y_i$
        for $1\le i<j$, $x_j'=z_{k-j}$ and $x_i'=y_{i-1}$ for $j<i\le k$. Thus the proof of \eqref{7.11-2} is complete.
        
        In particular, if $k>n$ then \eqref{7.11-2} applied to $j=n$ yields
        \[x_1\circ\dotsb\circ x_k=x_1'\circ\rho(x_2')\circ\dotsb\circ\rho^{n-1}(x_n')\circ\rho^n(x_{n+1}')\circ\dotsb\circ\rho^n(x_k').\]
        This leads to $K[M]=\sum_{f\in F}f\circ R$, where \[F=\{x_1\circ\dotsb\circ x_k:x_1,\dotsc,x_k\in X\text{ and }0\le k\le n\}\]
        is the subset of $M$ consisting of elements of degree $\le n$. So clearly $F$ is a finite set and the first part of the statement is proven.
    \end{proof}
\end{prop}

The class of finite solutions considered in Proposition~\ref{7.11} can be rewritten in terms of a generalization
of the well-studied algebraic system called a rack (see for instance \cite{MR3868941} and references therein).
Recall that a rack is a pair $(X,\tr)$, where the binary operation $\tr$ on the set  $X$ satisfies:
\begin{enumerate}
    \item the mapping $X\to X\colon y\mapsto x\tr y$ is a permutation for each $x\in X$,
    \item $x\tr(y\tr z)=(x\tr y)\tr(x\tr z)$ for all $x,y,z\in X$ (self-distributivity).
\end{enumerate}
It is well-known known \cite{MR3558231} that if $(X,\tr)$ is a rack then $(X,r)$, where $r(x,y)=(x\tr y,x)$,
is a bijective non-degenerate solution of the Yang--Baxter equation. More generally, consider triples  $(X,\tr,*)$,
where $X$ is a set and $X\times X\to X\colon(x,y)\mapsto x\tr y$ and $X\to X\colon x\mapsto x^*$ are mappings.
The triple is called a twisted rack if the following properties hold:
\begin{enumerate}
    \item the mapping $X\to X\colon y\mapsto x\tr y$  is a permutation of $X$ for each $x\in X$,
    \item $x\tr(y\tr z)=(x\tr y)\tr(x^*\tr z)$ for all $x,y,z\in X$ (twisted self-distributivity),
    \item $(x\tr y)^*=x^*\tr y^*$ for all $x,y\in X$.
\end{enumerate}
One easily verifies that  $(X,r)$, with $r(x,y)=(x\tr y,x^*)$, is a left non-degenerate solution of the Yang--Baxter equation. 
Moreover, each left non-degenerate solution $(X,r)$ of the form $r(x,y)=(\lambda_x(y),\rho(x))$ is determined by the twisted rack
$(X,\tr,*)$, where $x\tr y=\lambda_x(y)$ and $x^*=\rho(x)$. Hence the solutions considered in Proposition~\ref{7.11}
are in a one-to-one correspondence with finite twisted racks.
    
Obviously, if $G$ is a group and $f\colon G\to G$ is a metahomomorphism then defining $x\tr y=xyf(x)^{-1}$ and $x^*=f(x)$
one obtains a twisted rack $(G,\tr,*)$, which determines the solution $r(x,y)=(xyf(x)^{-1},f(x))$ mentioned in the introduction
of this section. Twisted conjugations $y\mapsto xyf(x)^{-1}$, in case $f$ is an automorphism of $G$, are being studied in
group theory, see for example \cite[Section 1.2]{MR2639839} or \cite{MR4434877,MR4243804}. 

\section*{Acknowledgement}

The third author is supported by National Science Centre grant 2020/39/D/ST1/01852 (Poland).
The fourth author is supported by Fonds voor Wetenschappelijk Onderzoek (Flanders), via an FWO post-doctoral fellowship, grant 12ZG221N.

\bibliographystyle{amsplain}
\bibliography{refs}
\end{document}